%Paper: q-alg/9505004
%From: pbueno@evalvx.ific.uv.es (J.CARLOS PIREZ)
%Date: Tue, 2 May 1995 19:04:41 +0200
%Date (revised): Wed, 3 May 1995 15:24:29 +0200
%Date (revised): Tue, 26 Sep 1995 08:39:46 +0200

%%%%%%%%%%%%%%%%%%%%%%%%%%%%%%%%%%%%%%%%%%%%%%%%%%%%%%%%%%%%%%%%%
%Version adaptada a las correcciones de estilo pedidas por el JMP:
%Accepted in JMP
%%%%%%%%%%%%%%%%%%%%%%%%%%%%%%%%%%%%%%%%%%%%%%%%%%%%%%%%%%%%%%%%%
\input phyzzx.tex
%%%%%%%%%%%%%%%%%%%%%%%%%%%%%%%%%%%%%%%%%%%%%%%%%%%%%%%%%%%%%%%%%%%
%  This is the PHYZZX.LOCAL file.
%  It contains local, ie. site-dependent macros.
%
%%%%%%%%%%%%%%%%%%%%%%%%%%%%%%%%%%%%%%%%%%%%%%%%%%%%%%%%%%%%%%%%%%%%
%
\catcode `\@=11
\newskip\frontpageskip
\newtoks\Ftuvnum   \let\ftuvnum=\Ftuvnum
\newtoks\Ificnum   \let\ificnum=\Ificnum
\newtoks\Pubtype  \let\pubtype=\Pubtype
\newif\ifp@bblock  \p@bblocktrue
\def\PH@SR@V{\doubl@true \baselineskip=18.1pt plus 0.2pt minus 0.1pt
             \parskip= 3pt plus 2pt minus 1pt }
\def\PHYSREV{\papers\PhysRevtrue\PH@SR@V}
\let\physrev=\PHYSREV
\def\titlepage{\FRONTPAGE\papers\ifPhysRev\PH@SR@V\fi
   \ifp@bblock\p@bblock \else\hrule height\z@ \rel@x \fi }
\def\nopubblock{\p@bblockfalse}
\def\endpage{\vfil\break}
\frontpageskip=12pt plus .5fil minus 2pt
\Pubtype={}
\Ftuvnum={}
\Ificnum={}
\def\p@bblock{\begingroup \tabskip=\hsize minus \hsize
   \baselineskip=1.5\ht\strutbox \topspace-2\baselineskip
   \halign to\hsize{\strut ##\hfil\tabskip=0pt\crcr
       \the\Ftuvnum\crcr\the\Ificnum\crcr
	\the\date\crcr\the\pubtype\crcr}\endgroup}

\newcount\anni
\anni=\year
\advance\anni by -1900
\def\aapub{\afterassignment\aap@b\toks@}
\def\aap@b{\edef\n@xt{\Ftuvnum={FTUV\ \the\anni--\the\toks@}}\n@xt}

\def\bbpub{\afterassignment\bbp@b\toks@}
\def\bbp@b{\edef\n@xt{\Ificnum={IFIC\ \the\anni--\the\toks@}}\n@xt}

\let\ftuvnum=\aapub
\let\ificnum=\bbpub
\def\ftuvbin{
      \expandafter\ifx\csname binno\endcsname\relax%
         \expandafter\ifx\csname MailStop\endcsname\relax%
         \else%
            , Mail Stop \MailStop%
         \fi%
      \else%
         , Mail Stop \binno%
      \fi%
   }
\expandafter\ifx\csname eightrm\endcsname\relax
    \let\eightrm=\ninerm  \fi
\def\memohead{\hrule height\z@ \kern -0.5in
    \line{\quad\fourteenrm FTUV MEMORANDUM\hfil \twelverm\the\date\quad}}
\def\memorule{\par \medskip \hrule height 0.5pt \kern 1.5pt
   \hrule height 0.5pt \medskip}
\def\FromAddress={
    Dr Moliner 50\cr Burjassot( Valencia)\cr }

\def\FTUVHEAD{\setbox0=\vtop{\baselineskip=10pt
     \ialign{\eightrm ##\hfil\cr
        Juan Carlos P\'{e}rez Bueno\cr
	Dr Moliner 50\cr
	Burjassot (Valencia)\cr
        C\'{o}digo Postal 46100\cr
        \cropen{1\jot}
        \ftuvphone\cr }}
   \setbox2=\hbox{\caps IFIC, Centro Mixto UV-CSIC}%
   \hrule height \z@ \kern -0.5in
   \vbox to 0pt{\vss\centerline{\fourteenrm DEPARTAMENTO DE F\'{I}SICA
TE\'{O}RICA}}
   \vbox{} \medskip
   \line{\hbox to 0.7\hsize{\hss \lower 10pt \box2 \hfill }\hfil
         \hbox to 0.25\hsize{\box0 \hfil }}\medskip }

\def\NAME{\author{J.A. de Azc\'{a}rraga \foot{azcarrag@evalvx.ific.uv.es}
and J.C. P\'{e}rez Bueno\foot{pbueno@evalvx.ific.uv.es}}}

\def\FTUV{\address{Departamento de F\'{\i}sica Te\'{o}rica and IFIC,\break
      Centro Mixto Univ. de Valencia-CSIC\break
	46100-Burjassot (Valencia), Spain}}

\def\ftuvphone{(96) 386--\ftuvext}
\def\ftuvext{4550}
\VOFFSET=33pt
\papersize
%
%%%%%%%%%%%%%%%%%%%%%%%%%%%%%%%%%%%%%%%%%%%%%%%%%%%%%%%%%%%%%%%%
%  Now comes the graphic package.
%  This version is rather primitive
%
\newwrite\figscalewrite
\newif\iffigscaleopen
\newif\ifgrayscale
\newif\ifreadyfile

\def\parsefilename{\ifreadyfile \else
    \iffigscaleopen \else \gl@bal\figscaleopentrue
       \immediate\openout\figscalewrite=\jobname.scalecon \fi
    \toks0={ }\immediate\write\figscalewrite{%
       \the\p@cwd \the\toks0 \the\p@cht \the\toks0 \the\picfilename }%
    \expandafter\p@rse \the\picfilename..\endp@rse \fi }
\def\p@rse#1.#2.#3\endp@rse{%
   \if*#3*\dop@rse #1.1..\else \if.#3\dop@rse #1.1..\else
                                \dop@rse #1.#3\fi \fi
   \expandafter\picfilename\expandafter{\n@xt}}
\def\dop@rse#1.#2..{\count255=#2 \ifnum\count255<1 \count255=1 \fi
   \ifnum\count255<10  \edef\n@xt{#1.PICT00\the\count255}\else
   \ifnum\count255<100 \edef\n@xt{#1.PICT0\the\count255}\else
                       \edef\n@xt{#1.PICT\the\count255}\fi\fi }
\def\redopicturebox{\edef\picturedefinition{\ifgrayscale
     \special{insert(\the\picfilename)}\else
     \special{mergeug(\the\picfilename)}\fi }}
%
%%%%%%%%%%%%%%%%%%%%%%%%%%%%%%%%%%%%%%%%%%%%%%%%%%%%%%
% Few miscellaneous macros
%

%
\def\pri{^{\, \prime }}

\def\boxit#1{\vbox{\hrule\hbox{\vrule\kern3pt
\vbox{\kern3pt#1\kern3pt}\kern3pt\vrule}\hrule}}
\catcode `\@=\active
%%%%%%%%%%%%%%%%%%%%%%%%%%%%%%%%%%%%%%%%%%%%%%%%%%%%%%%%%

\def\JMP#1(#2){\journal J. Math. Phys. &#1(#2)}
\def\PHL#1(#2){\journal Phys. Lett. &#1(#2)}
\def\JA#1(#2){\journal J. Algebra &#1(#2)}
\def\LMP#1(#2){\journal Lett. Math. Phys. &#1(#2)}
\def\JPH#1(#2){\journal J. Phys. &#1(#2)}

\REF\DRI{V. G. Drinfel'd, in Proc. of the 1986 {\it Int. Congr. of Math.}, MSRI
Berkeley, vol {\bf I}, 798 (1987) (A. Gleason, ed.)}

\REF\JIM{M. Jimbo {\it Lett. Math. Phys} {\bf 10}, 63 (1985);
{\it ibid.} {\bf 11}, 247 (1986)}

\REF\FRT{L.D. Faddeev, N. Yu. Reshetikhin and L. A. Takhtajan, {\it Alg. i
Anal.} {\bf 1}, 178 (1989) (Leningrad Math. J. {\bf 1}, 193 (1990))}

\REF\OSWZ{O. Ogievetsky, W. B. Schmidke, J. Wess and B. Zumino, {\it Commun.
Math. Phys.} {\bf 150}, 495 (1992)}

\REF\PODWOR{P. Podle\'{s} and S.L. Woronowicz, {\it On the classification of
quantum Poincar\'{e} groups}, hep-th/9412059}

\REF\WORZAK{S. L. Woronowicz and S. Zakrzewski, {\it Comp. Math.}
 {\bf 90}, 211 (1994)}

\REF\AKR{J. A. de Azc\'{a}rraga, P. P. Kulish and F. Rodenas
{\it Lett. Math. Phys} {\bf 32}, 173 (1994)}

\REF\AZRO{J. A. de Azc\'{a}rraga and F. Rodenas, {\it Deformed Minkowski
spaces: classification and properties}, Valencia preprint (August 1995)}

\REF\BHOS{A. Ballesteros, F. J. Herranz, M. A. del Olmo and M. Santander
{\it J. Math. Phys.}, {\bf 35}, 4928 (1994); {\it J. Math. Phys}
{\bf 27}, 1283 (1994)}

\REF\LNRT {J. Lukierski, A. Nowicki, H. Ruegg and V.N. Tolstoy,
{\it Phys. Lett.} {\it B264}, 331 (1991). A review is given in
J. Lukierski, H. Ruegg, and V.N. Tolstoy,
{\it $\kappa$-Quantum Poincar\'{e} 1994}, in {\it Quantum groups:
formalism and applications}, J. Lukierski, Z. Popowicz
and J. Sobczyk eds., PWN (1994), p. 359}

\REF\FIRI{E. Celeghini, R. Giachetti, E. Sorace and M. Tarlini
{\it J. Math. Phys.} {\bf 31}, 2548 (1990); {\it ibid.}
 {\bf 32}, 1155, 1159 (1991)}

\REF\MR{S. Majid and H. Ruegg {\it Phys. Lett.} {\bf B334}, 348 (1994)}

\REF\LODZ{S. Giller, P. Kosi\'nski, M. Majewski,
P. Ma\'slanka and J. Kunz {\it Phys. Lett.} {\bf B286}, 57 (1992) }

\REF\FOOTLUK{We thank J. Lukierski for pointing out Ref. \refmark{\LODZ}
to us.}

\REF\SIT{A. Sitarz {\it Phys. Lett.} {\bf B349}, 42 (1995)}

\REF\CONLOT{A. Connes and J. Lott, {\it Nucl. Phys} (Proc. Suppl.) {\bf B18},
29 (1990)}

\REF\KAST{D. Kastler and T. Sch\"ucker, {\it A detailed account of A. Connes
version of the standard model IV}, CPT-94/P.3092 (Jan 1995) hep-th 9501077}

\REF\MB{S. Majid {\it J. Algebra} {\bf 130}, 17 (1990)}

\REF\MAJSOB{S. Majid, {\it Israel J. Math.} {\bf 72}, 132 (1990);
S. Majid and Ya. S. Soibelman {\it J. Algebra} {\bf 163}, 68 (1994)}

\REF\FOOTDIMEN{The generators of a simple Lie algebra are dimensionless.
Although they may be
given dimensions by introducing constants in the r.h.s. of the commutators, a
re-scaling involving the constants makes them dimensionless (these constants
play the r\^ole of curvatures, and the re-scaling implies taking them as a
unit; see \refmark{\BHOSA} below).
The contraction, which `abelianizes' the algebra,
removes accordingly an internal scale and introduces the need of an
{\it external} dimension. Nevertheless, in a deformed \pk-theory the unit of
length $1/\kappa$ comes from $q$ and is therefore included in the algebra.}

\REF\BHOSA{A. Ballesteros, F. J. Herranz, M. A. del Olmo and M. Santander
{\it J. Phys.} {\bf A23}, 5801 (1993) }

\REF\AG{J. A. de Azc\'arraga and D. Ginestar {\it J. Math. Phys.}
{\bf 32}, 3500 (1991)}

\REF\KLMNS{J. Lukierski, J. Sobczyk and A. Nowicki {\it J. Phys.}
{\bf A26}, L1099 (1993);
P. Kosi\'nski, J. Lukierski, P. M\'aslanska and J. Sobczyk {\it J. Phys.}
{\bf A28}, 2255 (1995)}

\REF\CHANG{Zhe Chang, Wei Cheng and Han-Ying Guo, {\it J. Phys.} {\bf A23},
4185 (1990)}

\REF\FOOTQGR{The name `quantum group' is borrowed from quantum physics due
to the noncommutative geometry implied by the deformation. Since it is
convenient to distinguish among both effects we shall say  \eg,
`deformed' or `noncommutative' rather than `quantum' spacetime.}

\REF\AKRA{J. A. de Azc\'arraga, P. P. Kulish and F. Rodenas,
{\it Phys. Lett.} {\bf B351}, 23 (1995)}

\REF\FOOTBRAID{Majid's braided geometry corresponds
to this second approach. For an introduction, see S. Majid,
{\it Introduction to braided geometry and q-Minkowski
space}, Proc. of the Varenna 1994 School in quantum groups, hep-th 9410241}

\REF\FOOTTRANS{We note in passing that ${\cal U}_\kappa(Tr)$ has itself a
bicrossproduct structure,
${\cal U}_\kappa(Tr)={\cal U}(Tr_3)\bic{\cal U}(Tr_0)$, where the two Hopf
algebras are generated by the abelian generators $P_0$ and $P_i$
respectively with primitive coproducts, $\alpha$ is
trivial ($P_0\acti P_i=[P_0,P_i]=0$) and
$\beta(P_i)=\exp(-P_0/\kappa)\otimes P_i$.
Specifically, ${\cal U}_\kappa(Tr)$ is a semidirect coproduct of Hopf
algebras (Appendix A).}

\REF\LUKRUE{J. Lukierski and H. Ruegg {\it Phys. Lett.} {\bf B329}, 189
(1994)}

\REF\SZAK{S. Zakrzewski {\it J. Phys.} {\bf A27}, 2075 (1994)}

\REF\FOOTLEIB{This avoids
the left exponential factor $\exp(-P_0/\kappa)$ in \tri.
In general one has \eg,
$P_i\act(ff\pri):=m(\Delta(P_i)\act(f\otimes f\pri))\equiv(P_i\act f)
f\pri+(\exp(-P_0/\kappa)\act f)(P_i\act f\pri)$.
Thus, and due to its non-primitive coproduct, $P_i$ may be identified
with $\partial_i$ and will satisfy Leibniz's rule only if the above
condition is met: if $f\neq f(x_0),\ P_i\act(ff\pri)=(\partial_i f)f\pri
+f(\partial_i f\pri)$.}

\REF\WOR{S. L Woronowicz, {\it Commun. Math. Phys.} {\bf 122}, 125 (1989)}

\REF\FOOTPHI{The $\phi$ in \difxi\ differs from that in
\refmark{\SIT} by a
factor \kap\ so that in the limit $\kappa\rightarrow\infty$
 all relations \difxi\ become commutative.  Thus, our
one-form $\phi$ has dimensions of length.}

\REF\AACON{V. Aldaya and J.A. de Azc\'{a}rraga, {\it Int. J. of Theor. Phys.}
{\bf 24}, 141 (1985)}

\REF\AI{J. A. de Azc\'{a}rraga and J. M. Izquierdo, {\it Lie algebras, Lie
groups cohomology and some applications in physics}, Camb. Univ. Press,
to appear (1995)}

\REF\BAR{V. Bargmann, {\it Ann. Math.} {\bf 59}, 1 (1954)}

\REF\SAL{E. J. Saletan {\it J. Math. Phys} {\bf 2}, 1 (61)}

\REF\FOOTCOCH{The one-cochain generating the two-coboundary
defining ${\cal P}\times
U(1)$ corresponds to the subtraction of the rest energy ($mc^2$) from $p^0c$
and to the redefinition of the Klein-Gordon wavefunction,
steps which are needed before performing the nonrelativistic limit.}

\REF\IW{E. \.In\"on\"u and E.P. Wigner, {\it Proc. Nat. Acad. Sci.} {\bf 39},
510 (1953)}

\REF\FOOTFACTOR{Note, however, that in
contrast with \refmark{\LODZ} we are {\it not} considering
the $c\rightarrow\infty$ limit of \pk\ for $\kappa=\tilde\kappa/c$,
$\tilde\kappa$ constant.}

\REF\CWSWW{U. Carow-Watamura, M. Schlieker, S. Watamura and W. Weich, {\it
Commun. Math. Phys.} {\bf 142}, 605 (1991)}

\REF\ZSAL{B. Zumino {\it Differential calculus on quantum spaces and quantum
groups}, in XIX ICGTMP, Salamanca (Anales de F\'{\i}sica (monograf\'{\i}as)
{\bf 1}, vol. I, p. 44 (1992))}

\REF\SCH{K. Schm\"udgen and A. Sch\"uler, {\it C.R. Acad. Sci. Paris}
{\bf 316}, 1155 (1993)}

\REF\AS{A. Sudbery, {\it The quantum orthogonal mystery} in {\it Quantum
groups} (Karpacz 1994), J. Lukierski, Z. Popowicz and J. Sobczyk eds., PWN
(1995), p. 303}

\REF\BCM{R. J. Blattner, M. Cohen and S. Montgomery, {\it Trans. Am. Math.
Society} {\bf 298}, 671 (1986); R. J. Blattner and S. Montgomery, {\it Pac. J.
Math.} {\bf 137}, 37 (1989)}

\REF\FOOTCOND{If $H$ is cocommutative and $A$ commutative,
condition \apxiii\ is automatically satisfied. This is always the
case in the main text, where $A$ corresponds to translations and $H$ to
rotations and boosts.}

\REF\MOL{R. Molnar {\it J. Algebra} {\bf 47}, 29 (77)}

\REF\FOOTASSOC{This follows from associativity. For the Hopf algebra
with primitive coproduct defined on the enveloping algebra
${\cal U}({\cal G})$
of a Lie algebra ${\cal G}$, eq. \apxxi\ gives
the familiar two-cocycle condition for the right action of $H$ on $A$,
$\xi([h,g],f)-\xi([h,f], g)+\xi([g,f],h)+
\xi(h, g)\acti f- \xi(h,f)\acti
g+\xi(g,f)\acti h=0,\;h,g,f\in{\cal G}.$}

%\REF\INONU{E. \.In\"on\"u, {\it Contractions of Lie groups and their
%representations} in {\it Group theor. concepts in elem. part. physics}, F.
%G\"ursey ed., Gordon and Beach, p. 391 (1964)}
%%%%%%%%%%%%%%%%%%%%%%%%%%%%%%%%%%%%%%%%%%%%%%%%%%
\catcode `\@=11
\def \poiv {({\Number {2}}.1)}
\def \pov {({\Number {2}}.2)}
\def \poiii {({\Number {2}}.3)}
\def \tri {({\Number {2}}.4)}
\def \aldef {({\Number {2}}.5)}
\def \bedef {({\Number {2}}.6)}
\def \IIxx {({\Number {2}}.7)}
\def \IIxxi {({\Number {2}}.8)}
\def \difi {({\Number {3}}.1)}
\def \difii {({\Number {3}}.2)}
\def \dconsis {({\Number {3}}.3)}
\def \difiii {({\Number {3}}.4)}
\def \difv {({\Number {3}}.5)}
\def \difvi {({\Number {3}}.6)}
\def \difvii {({\Number {3}}.7)}
\def \difviii {({\Number {3}}.8)}
\def \difix {({\Number {3}}.9)}
\def \difa {({\Number {3}}.10)}
\def \difxi {({\Number {3}}.11)}
\def \NUU {({\Number {3}}.12)}
\def \difxii {({\Number {3}}.13)}
\def \difxiii {({\Number {3}}.14)}
\def \difxiv {({\Number {3}}.15)}
\def \difxv {({\Number {3}}.16)}
\def \VIIi {({\Number {4}}.1)}
\def \VIIa {({\Number {4}}.2)}
\def \VIIb {({\Number {4}}.3)}
\def \VIIii {({\Number {4}}.4)}
\def \VIIiii {({\Number {4}}.5)}
\def \VIIiv {({\Number {4}}.6)}
\def \VIIv {({\Number {4}}.7)}
\def \VIIva {({\Number {4}}.8)}
\def \VIIvi {({\Number {4}}.9)}
\def \rescal {({\Number {5}}.1)}
\def \Bi {({\Number {5}}.2)}
\def \Bii {({\Number {5}}.3)}
\def \antipod {({\Number {5}}.4)}
\def \extras {({\Number {5}}.5)}
\def \Biii {({\Number {5}}.6)}
\def \Biv {({\Number {5}}.7)}
\def \Bv {({\Number {5}}.8)}
\def \Bvi {({\Number {5}}.9)}
\def \Bvii {({\Number {5}}.10)}
\def \Bix {({\Number {5}}.11)}
\def \Bx {({\Number {5}}.12)}
\def \gali {({\Number {6}}.1)}
\def \galii {({\Number {6}}.2)}
\def \galiii {({\Number {6}}.3)}
\def \galiv {({\Number {6}}.4)}
\def \galv {({\Number {6}}.5)}
\def \galvi {({\Number {6}}.6)}
\def \defalfa {({\Number {6}}.7)}
\def \defbeta {({\Number {6}}.8)}
\def \demiii {({\Number {6}}.9)}
\def \casi {({\Number {6}}.10)}
\def \casii {({\Number {6}}.11)}
\def \casiii {({\Number {6}}.12)}
\def \vi {({\Number {7}}.1)}
\def \viii {({\Number {7}}.2)}
\def \viv {({\Number {7}}.3)}
\def \vv {({\Number {7}}.4)}
\def \vvi {({\Number {7}}.5)}
\def \api {({A}.1)}
\def \apii {({A}.2)}
\def \apiiia {({A}.3)}
\def \apiii {({A}.4)}
\def \apiv {({A}.5)}
\def \apv {({A}.6)}
\def \apvi {({A}.7)}
\def \apvii {({A}.8)}
\def \apviii {({A}.9)}
\def \apix {({A}.10)}
\def \apx {({A}.11)}
\def \apxi {({A}.12)}
\def \apxii {({A}.13)}
\def \apxiii {({A}.14)}
\def \apxiv {({A}.15)}
\def \apxv {({A}.16)}
\def \apxvi {({A}.17)}
\def \apxvii {({A}.18)}
\def \apxx {({A}.19)}
\def \apxxi {({A}.20)}
\def \apxxii {({A}.21)}
\def \apxxiii {({A}.22)}
\def \apxxiv {({A}.23)}
\def \apxxv {({A}.24)}
\def \apxxvi {({A}.25)}
\def \apxxvii {({A}.26)}
\def \apxxviii {({A}.27)}
\def \apxxviiia {({A}.28)}
\def \apxxix {({A}.29)}
\def \apxxx {({A}.30)}
\catcode`\@=\active
%%%%%%%%%%%%%%%%%%%%%%%%%%%%%%%%%%%%%%%%%%%%%%%%%%
\catcode`\@=11

\font\tenmsa=msam10
\font\sevenmsa=msam7
\font\fivemsa=msam5
\font\tenmsb=msbm10
\font\sevenmsb=msbm7
\font\fivemsb=msbm5
\newfam\msafam
\newfam\msbfam
\textfont\msafam=\tenmsa  \scriptfont\msafam=\sevenmsa
  \scriptscriptfont\msafam=\fivemsa
\textfont\msbfam=\tenmsb  \scriptfont\msbfam=\sevenmsb
  \scriptscriptfont\msbfam=\fivemsb

\def\hexnumber@#1{\ifnum#1<10 \number#1\else
 \ifnum#1=10 A\else\ifnum#1=11 B\else\ifnum#1=12 C\else
 \ifnum#1=13 D\else\ifnum#1=14 E\else\ifnum#1=15 F\fi\fi\fi\fi\fi\fi\fi}

\def\msa@{\hexnumber@\msafam}
\def\msb@{\hexnumber@\msbfam}
\mathchardef\blacktriangleright="3\msa@49
\mathchardef\blacktriangleleft="3\msa@4A
\catcode`\@=\active

\newlist\mypoint=\alphabetic&)&0.5\itemsize;
\def\pk{${\cal P}_\kappa$}
\def\kap{$\kappa$}
\def\kapp{$\tilde\kappa$}

\def\bic{\triangleright\!\!\!\blacktriangleleft}
\def\acti{\triangleleft}
\def\act{\triangleright}
\def\rimo{\triangleright\!\!\!<}
\def\leco{>\!\!\blacktriangleleft}
\def\gk{${\cal G}_{\tilde\kappa}$}
\def\gkk{$\tilde{\cal G}_{(m)\hat\kappa}$}
\def\nkk{${\cal N}_{\hat\kappa}$}
\def\and{\times}
%%%%%%%%%%%%%%%%%%%%%%%%%%%%%%%%%%%%%%%%%%%%%%%%%%%%%
\physrev

\ftuvnum={12}
\ificnum={12}
\pubtype={q-alg/9505004}
\date={April 20, 1995}
\titlepage
\hfill To appear in JMP

\title{Relativistic and Newtonian \kap-spacetimes}\foot{ PACS:
02.20.+b, 02.90.+p, 03.30.+p, 11.30.Cp }
\NAME
\FTUV
\abstract
 The deformations of the Galilei algebra and their associated noncommutative
 Newtonian spacetimes are investigated. This is done by analyzing the possible
 nonrelativistic limits of an eleven generator (pseudo)extended
\kap-Poincar\'e algebra $\tilde{\cal P}_\kappa$
 and their implications for the existence of a first order differential
 calculus. The additional one-form needed to achieve a consistent calculus on
 \kap-Minkowski space is shown to be related to the additional central
 generator entering in the $\tilde{\cal P}_\kappa$ Hopf algebra.
 In the process, deformations of the
 extended Galilei and Galilei algebras are introduced which have,
 respectively, a cocycle and a bicrossproduct structure.\endpage
\chapter{Introduction}
The deformation or `quantization' of the symmetry Lie group of affine spaces
has been hindered by the lack of a prescription like the Drinfel'd-Jimbo one
\refmark{\DRI-\JIM} which applies to simple groups and for which there
 is also a well defined universal $R$-matrix \refmark{\FRT}.
 From a physical point of view, the most interesting groups to deform are
 the kinematical groups of relativistic and nonrelativistic theories, the
 Poincar\'{e} and the Galilei groups. But due to the lack of a definite
 prescription for inhomogenous groups, there is not a unique deformed
 Poincar\'{e} algebra. A recent classification of deformed Poincar\'e
groups (which nevertheless does not include all proposals as \eg,
\refmark{\OSWZ}) has been given in
\refmark{\PODWOR} based in the deformations of the Lorentz group
\refmark{\WORZAK} (see also \refmark{\AKR}\refmark{\AZRO}
in connection with
deformed Minkowski spaces). Other deformed spacetime affine
algebras have been proposed in \refmark{\BHOS}.
We shall devote this paper to the problem of defining a deformation
of the Galilei algebra and its associated spacetime. Our starting point will be
the \kap-Poincar\'{e} algebra \pk\ \refmark{\LNRT},
 which is obtained by a non-standard contraction
 (\ie, involving also the
 deformation parameter $q$ \refmark{\FIRI}) from the deformed
anti-De Sitter algebra ${\cal U}_q(so(3,2))$.
{}From a mathematical point of view, \pk\  has the interest of having
a bicrossproduct structure
 \refmark{\MR}, which is specially adept to deform Lie groups with a
semidirect character and hence inhomogeneus kinematical groups; it is also one
of the deformations in \refmark{\BHOS} and
\refmark{\PODWOR}. We shall obtain deformed
Galilean algebras and Newtonian spacetimes by analyzing the contractions of
\pk. The second Galilean deformation uncovered by our analysis, denoted
\gk\ in sec. 6, was given in \refmark{\LODZ} in a different basis
\refmark{\FOOTLUK}. The bicrossproduct
structure and the Casimirs of \gk\ (sec. 6) were not,
however, discussed in this reference.

A bicovariant and Lorentz covariant first order differential calculus
on \kap-\break Minkowski spacetime ${\cal M}_\kappa$
has been proposed recently
 \refmark{\SIT}. The self-consistency of this
 differential calculus requires the addition of one
 scalar one-form $\phi$ to the spacetime ones $dx_\mu$. We show that this
additional variable may be related to a new central generator which
determines a Hopf algebra `pseudoextension' (see below)
of \pk. To stress the
analogy with the undeformed case, we study a differential calculus based
on a non-scalar form $\varphi$, which in the nonrelativistic
 limit leads to a consistent differential calculus on an enlarged
Newtonian spacetime associated with a deformation of the extended Galilei
algebra. Although, up to now, there seems to be no physical need for deforming
 spacetime and the applications of noncommutative geometry to real physical
 theories do not abound (see, however, \refmark{\CONLOT-\KAST}),
 the analysis of the above problems may shed some light on the nature of the
 deformation process.

All the algebras considered in this paper have a bicrossproduct \refmark{\MB}
or a cocycle bicrossproduct structure \refmark{\MAJSOB}. The defining
properties of these structures are summarized for completeness
in Appendix A.

\vskip 24 pt
\chapter{ ${\cal P}_\kappa$
and \kap-Minkowski spacetime}
Let us start by recalling the defining relations
 of the \kap-Poincar\'{e} algebra \pk
\refmark{\LNRT} in the basis which is
suitable to exhibit the bicrossproduct structure \refmark{\MR}
$$\eqalign{&[P_\mu,P_\nu]=0\quad,\quad [M_i,M_j]=\epsilon_{ijk}M_k\quad,\quad
[M_i,P_j]=\epsilon_{ijk}P_k\quad,\cr
&[M_i,P_0] =0\quad,\quad
[M_i, N_j] =\epsilon_{ijk}N_k\quad,\quad
[N_i,P_0] =P_i\quad,\cr
&[N_i,N_j] =-\epsilon_{ijk}M_k\quad,\cr
&[N_i,P_j] =\delta_{ij}\left[{\kappa\over 2} \left(1-\exp\left(-{2P_0
\over\kappa}\right)\right)+{1\over{2\kappa}}{\bf P}^2\right]-
{1\over\kappa} P_iP_j\quad;\cr}\eqn\poiv$$
this basis preserves the classical Lorentz subalgebra.
The Hopf algebra structure of \pk\ is given by adding the coproducts
$$\eqalign{\Delta P_0 =&P_0\otimes 1+1\otimes P_0\quad,\quad
\Delta  P_i =P_i\otimes 1+\exp\left({-P_0\over
\kappa}\right)\otimes P_i\quad, \cr
\Delta M_i =&M_i\otimes 1+1\otimes M_i \quad,\cr
\Delta N_i =&N_i\otimes 1+
\exp\left({-P_0\over\kappa}\right)\otimes N_i +{1\over\kappa} \epsilon_{ijk}
 P_j\otimes M_k\quad;\cr}\eqn\pov$$
counits ($\epsilon(P_\mu,M_i,N_j)=0$) and antipodes
$$\eqalign{&S(P_0) =-P_0\quad,\quad
S(P_i) =-\exp\left({P_0\over\kappa}\right)P_i\quad, \cr &
S(M_i) =-M_i\quad,\quad
S(N_i) =-\exp\left({P_0\over\kappa}\right)N_i +{1\over\kappa} \epsilon_{ijk}
\exp\left({P_0\over\kappa}\right)P_jM_k\quad. \cr}\eqn\poiii$$

The deformation parameter \kap, with dimensions of (length)$^{-1}$, appears
after contracting ${\cal U}_q(so(2,3))$ by rescaling
$M_{5\mu}=RP_{\mu}$
{\it and}
 redefining the original deformation parameter $q$ as $q=\exp(1/\kappa R)$
before performing the $R\rightarrow\infty$ limit; thus
 $M_i$ and $N_i$ are dimensionless and $[P_\mu]=$ (length)$^{-1}$.
 As a result, \pk\ introduces a `natural' {\it length} unit $1/\kappa$
in a \pk-governed relativistic
theory \refmark{\FOOTDIMEN}.
Since \kap\ appears as a new `universal' constant still undetermined,
it is possible to think of \kap\ as including factors of the other natural
constant in the theory, the velocity of light $c$. In order to discuss
possible nonrelativistic limits, we shall consider below two possibilities:
a) replacing $\kappa$ by $\hat\kappa c$ and
b) replacing $\kappa$ by ${\tilde\kappa}/c$,
which correspond
to deformation parameters with dimensions $[\hat\kappa]=L^{-2}T$,
$[\tilde\kappa]={\rm frequency}$. Notice that we do {\it not} set
$\kappa=\hat\kappa c$ or $\kappa=\tilde\kappa/c$ in \pk; rather, we shall
consider deformations ${\cal P}_{\hat\kappa}$ and ${\cal P}_{\tilde\kappa}$ of
${\cal P}$ {\it defined} by making in \pk\ the above replacements.
In fact, any factors in $c$ hidden by the use of units in which $c=1$ have
to be made explicit to discuss the non-relativistic limit. They may appear
accompanying constants such as $\kappa$ here or in other places, \eg\ as in
the cocycle defining the supertranslation graded group where an $1/c$
factor is needed to define the non-relativistic limit of supersymmetry
\refmark{\AG}.
This fact may be used to obtain, from the bicrossproduct structure of the
$\kappa$ deformed superPoincar\'e algebra \refmark{\KLMNS}, the corresponding
deformed superGalilei algebra.

We do not use natural units and Planck's constant $\hbar$
will not appear in the text. Thus, $[\kappa]=L^{-1}$; \kap\
cannot have dimensions of mass in a classical ($\hbar=0$) framework.
There are two general approaches to discuss deformation ($q\not =1$) and
quantization ($\hbar\not = 1$) (see \refmark{\CHANG} for an early discussion
of both processes).
If $q=q(\hbar)$ is assumed (and there is no a priori reason for it)
mathematical deformation (also termed `quantization' \refmark{\FOOTQGR})
implies physical quantization. Since $q$ is dimensionless, this requires
the presence of another dimensional constant that will survive in the
quasiclassical approximation which itself requires
a {\it definite} hypothesis on the form of $q(\hbar)$
(as \eg, $q=\exp(\gamma\hbar)$). It constitutes an interesting problem to look
in this approximation at the effects of $\gamma$ in spacetime theories
(for a simple analysis already exhibiting the serious difficulties that may be
encountered see \refmark{\AKRA}). If $q\not =q(\hbar)$, quantization and
deformation are different processes, and the presence of $q\not =1$ affects
the commuting properties of the algebra elements at the {\it
classical} level and deforms the `product' between the algebras
\refmark{\FOOTBRAID}.
There is no clear way to quantize (\ie, to introduce $\hbar$) a
deformed system, and the heuristic multiplication of the classical generators
by $i\hbar$, which can be justified geometrically in the undeformed Lie
algebra case,
might not be an adequate prescription to obtain  quantum operators
when $q\not =1$. We shall take here the $q\not =q(\hbar)$
point of view and restrict
ourselves essentially to classical considerations. Thus, the translation
generators will have dimensions of $L^{-1}$ rather than of momenta, and the
mass dimension will be introduced through the parameter characterizing
the two-cocycle of a central extension.

To recover the Hopf algebra \poiv-\poiii\
  as the bicrossproduct \refmark{\MB}
${\cal U}(so(1,3))\bic {\cal U}_\kappa(Tr)$ \refmark{\MR}
of Hopf algebras, one needs a right action
$\alpha:{\cal U}_\kappa (Tr)\otimes {\cal U}(so(1,3))\rightarrow
{\cal U}_\kappa (Tr)$,
characterizing ${\cal U}_\kappa (Tr)$ as a right $so(1,3)$-module algebra and
a left coaction
$\beta:{\cal U}(so(1,3))\rightarrow {\cal U}_\kappa (Tr)\otimes
{\cal U}(so(1,3))$,
characterizing ${\cal U}(so(1,3))$ as a left ${\cal U}_\kappa(Tr)$-comodule
coalgebra,
subjected to certain compatibility conditions (see Appendix A).
The algebra ${\cal U}_\kappa (Tr)$ is the deformation of
the translation Hopf algebra defined by:
$$\eqalign{& [P_\mu,P_\nu]=0\quad,\quad
 \Delta P_0=P_0\otimes 1+1\otimes P_0\quad, \cr &
 \Delta P_i=P_i\otimes 1 + \exp(-P_0/\kappa)\otimes P_i\quad; \cr &}\eqn\tri$$
${\cal U}(so(1,3))$ is the cocommutative Hopf algebra defined
on the enveloping algebra of the Lorentz algebra.
The right action of the {\it classical} Lorentz algebra on $Tr$,
$\alpha(t,h)\equiv t\acti h,\ t\in Tr,
 h\in\cal L$ is defined to be (cf. \poiv)
$$\eqalign{&P_0\acti M_i \equiv[P_0,M_i]=0\quad,\quad
P_i\acti M_j \equiv[P_i,M_j]=\epsilon_{ijk}P_k\quad,\cr &
P_0\acti N_i \equiv[P_0,N_i]=-P_i\quad,\cr &
P_i\acti N_j \equiv[P_i,N_j]=
-\delta_{ij}\left[{\kappa\over 2}(1-\exp(-2P_0/\kappa))
+{1\over{2\kappa}}{\bf P}^2\right]+{1\over\kappa} P_iP_j\quad, \cr
&}\eqn\aldef$$
and the left coaction $\beta$ is given by
$$\eqalign{&\beta(M_i) =1\otimes M_i\quad,\quad
\beta(N_i) =\exp(-P_0/\kappa)\otimes N_i + {\epsilon_{ijk}\over\kappa}
P_j\otimes M_k\quad. \cr &}\eqn\bedef$$
Using \apxv, \apxvii\ it is seen that \pov, \poiii\ are recovered and hence
\pk $={\cal U}(so(1,3))\bic {\cal
U}_\kappa(Tr)$\refmark{\MR}\refmark{\FOOTTRANS}.

The associated {\it \kap-Minkowski spacetime}
algebra ${\cal M}_\kappa$ may now be
introduced \refmark{\MR} as the dual $Tr^*$ to the translation (momentum)
sector $Tr$,
$<\!P_\mu,x^\nu\!>=\delta_\mu^\nu$. The commutativity
(noncocommutativity) of ${\cal U}_\kappa(Tr)$ induces cocommutativity
(noncommutativity)
in ${\cal M}_\kappa$. Specifically \refmark{\LUKRUE,\SZAK,\MR}
$$\Delta x_\mu=x_\mu\otimes 1+1\otimes x_\mu\quad;\quad
[x_i,x_j]=0\quad,\quad [x_i,x_0]={x_i\over\kappa}\quad.\eqn\IIxx$$
The canonical action of the momenta generators on ${\cal M}_\kappa$ is now
defined by $t\act t^*=<t^*_{(1)},t>t^*_{(2)}$; this leads to
$P_\mu=\partial/\partial x^{\mu}$ provided it acts on elements of
${\cal M}_\kappa$ with all powers of $x^0$ to the right \refmark{\MR}
\refmark{\FOOTLEIB}. As for the elements of
${\cal U}(so(1,3))$, their right action on $Tr$ induces a left one on
${\cal M}_\kappa$ by $<a\acti h,x>=<a,h\act x>$. Then eqs. \aldef\ lead to
$$M_i\act x_j=\epsilon_{ijk}x_k\quad,\quad M_i\act x_0=0\quad,\quad
N_i\act x_j=-\delta_{ij}x_0\quad,\quad N_i\act x_0=-x_i\eqn\IIxxi$$
(to which one may add $P_\mu\act x_\nu=\delta_{\mu\nu}$). Using
 these results on quadratic terms
($h\act xy=(h_{(1)}\act x)(h_{(2)}\act y)$) it is found that
$x^2_0-{\bf x}^2+{3\over\kappa}x_0$ is Lorentz invariant \refmark{\MR}
(\kap-Minkowski metric; however, it is not a central element in ${\cal
M}_{\kappa}$).

\chapter{Lorentz covariant differential calculus on ${\cal M}_{\kappa}$}
A feature of the Lorentz covariant calculus on
${\cal M}_\kappa$ of \refmark{\SIT}
is that the spacetime algebra needs to be enlarged with the addition
of a one-form $\phi$; otherwise there is no consistent
 solution for the relations defining
the bicovariant calculus. Given a Hopf algebra ${\cal A}$,
a first order bicovariant differential
calculus over ${\cal A}$ is defined \refmark{\WOR} by a pair $(\Gamma, d)$
where $d:{\cal A}\to\Gamma$ is a linear mapping satisfying Leibniz's rule and
$\Gamma$ is a bicovariant ${\cal A}$-bimodule \ie\
the linear mappings $\Delta_L:
\Gamma\rightarrow{\cal A}\otimes\Gamma\;,\;
\Delta_R:\Gamma\rightarrow\Gamma\otimes{\cal A}$ (left and right coactions)
and the exterior derivative $d$ satisfy
$$\eqalign{\Delta_L(a\omega b)=\Delta(a)\Delta_L(\omega)\Delta(b)\quad, &\quad
\Delta_R(a\omega b)=\Delta(a)\Delta_R(\omega)\Delta(b)\quad,\cr
(\Delta\otimes id)\Delta_L=(id\otimes\Delta_L)\Delta_L\quad,&\quad
(id\otimes\Delta)\Delta_R=(\Delta_R\otimes id)\Delta_R\quad,\cr}\eqn\difi$$
$$\eqalign{
(id\otimes\Delta_R)\Delta_L=&(\Delta_L\otimes id)\Delta_R\quad,\cr}\eqn\difii$$
$$\Delta_L d=(id\otimes d)\Delta\quad,\quad\Delta_R d=(d\otimes id)\Delta\quad,
\eqn\dconsis$$
where the left (right) equations in \difi\ express the left- (right-)
covariance of
$\Gamma$, \difii\ is the result of bicovariance (commutation of $\Delta_R$ and
$\Delta_L$) of $\Gamma$ and
\dconsis\ expresses the compatibility of $d$ and $\Delta$ with
$\Delta_L\;,\;\Delta_R$.
It follows from $\Delta x_\mu$ (eq. \IIxx) that all $dx_\mu$ are left-
(LI) and right-invariant \ie,
$$\Delta_L(dx_\mu)=1\otimes dx_\mu\quad,\quad
\Delta_R(dx_\mu)=dx_\mu\otimes 1\quad.\eqn\difiii$$
Following \refmark{\SIT}, let $\chi_a\quad (a=0,1,...,N\!\ge\! D-1,\ {\rm
where}\ D\ $is the dimension of spacetime) be a basis of left-invariant forms.
It then follows that the commutator $[x_\mu,\chi_a]$ is LI,
$\Delta_L([x_\mu,\chi_a])=1\otimes [x_\mu,\chi_a],$
 and hence that
$$[x_\mu,\chi_a]=A^{b}_{\mu a}\chi_{b}\quad.\eqn\difv$$
The Jacobi identity for $(x_\mu,x_\nu,\chi_a)$ then gives
$$B_{\mu\nu}^\rho A_{\rho a}^c+A_{\nu a}^bA_{\mu b}^c+A_{\mu a}^bA_{\nu b}^c=0
\quad,\eqn\difvi$$
where the commutators in \IIxx\ have been jointly expressed as
$[x_\mu,x_\nu]=B_{\mu \nu}^\rho x_\rho$. Since $dx_\mu$ is LI,
 $dx_\mu=C_\mu^a
\chi_a$, and Leibniz's rule applied to $[x_\mu,x_\nu]$ gives
$$-C_\mu^aA_{\nu a}^c+C_{\nu}^aA_{\mu a}^c=B_{\mu\nu}^\rho C_\rho^c\quad.
\eqn\difvii$$
The solutions to \difvi\ and \difvii\ determine a first order bicovariant
differential
calculus on ${\cal M}_\kappa$.

The action \IIxxi\ of the Lorentz algebra is extended to the module
of one-forms in the natural way
$$h\act(ydx)=(h_{(1)}\act y)(d(h_{(2)}\act x))\quad,\quad h\act (dxy)=
(d(h_{(1)}\act x))(h_{(2)}\act y)\quad.\eqn\difviii$$
This leads to the following relations
$$\eqalign{&N_k\act[x_i,dx_j]=-\delta_{ki}[x_0,dx_j]-\delta_{kj}[x_i,dx_0]+
{1 \over \kappa}(\delta_{kj}dx_i-\delta_{ij}dx_k)\quad,\cr
&N_k\act[x_0,dx_i]=-[x_k,dx_i]-\delta_{ki}[x_0,dx_0]+{1\over \kappa}\delta_{ki}
dx_0\quad,\cr
&N_k\act[x_i,dx_0]=-[x_i,dx_k]-\delta_{ki}[x_0,dx_0]\quad,\cr
&N_k\act[x_0,dx_0]=-[x_k,dx_0]-[x_0,dx_k]+{1 \over \kappa}dx_k\quad,\cr
&M_k\act[x_i,dx_j]=\epsilon_{kil}[x_l,dx_j]+\epsilon_{kjl}[x_i,dx_l]\quad,\quad
M_k\act[x_0,dx_i]=\epsilon_{kil}[x_0,dx_l]\quad,\cr
&M_k\act[x_i,dx_0]=\epsilon_{kil}[x_l,dx_0]\quad,\quad
M_k\act[x_0,dx_0]=0\quad,\cr}\eqn\difix$$
Now, since in general $h\act[x_\mu,\chi_a]=A_{\mu a}^c (h\act\chi_c)$ by eq.
\difv, we may set $\chi_\mu=dx_\mu$ and look for solutions to the system of
equations to which eqs. \difix\ give rise.
In fact, the solution
$$[x_i,dx_j]=\delta_{ij}{dx_0\over\kappa}\quad,\quad
[x_i,dx_0]={dx_i\over\kappa}\quad,\quad
[x_0,dx_i]=0\quad,\quad
[x_0,dx_0]=0\eqn\difa$$
is unique, but since it does not satisfy \difvi, it does
not define a covariant calculus \refmark{\SIT}. To obtain a consistent
solution,
an additional scalar ($M_i\act\phi=0=N_i\act\phi$) one-form
is necessary \refmark{\SIT},
 which leads to the solution \refmark{\FOOTPHI}

$$\eqalign{&[x_\mu,\phi]={1\over\kappa}dx_\mu\quad,
\quad [x_0,dx_0]={\phi\over\kappa}
\quad,\quad [x_0,dx_i]=0\quad,\cr
&[x_i,dx_0]={1\over\kappa}dx_i\quad,\quad
[x_i,dx_j]=\delta_{ij}{1\over\kappa}(dx_0-\phi)\cr}\eqn\difxi$$
which satisfies both \difvi\ and \difvii.

Two questions immediately arise. The first is the origin of the additional
one-form $\phi$; the
second is related with the possibility of defining a nonrelativistic limit
of \pk\ since eqs. \difa, \difxi\ do not have a $c\rightarrow\infty$ limit
$(x_0\equiv ct$) unless \kap\ is redefined. There are now two possibilities:

a) if we replace $\kappa$ by $\hat\kappa c$,
all commutators in \difa\ or \difxi\
become zero for $c\rightarrow\infty$ with the exception of
$$[x_i,dx_j]=\delta_{ij}{dt\over\hat\kappa}\quad;\eqn\NUU$$

b) we may replace $\phi$ by a different one-form $\varphi$ which, in
contrast with the scalar $\phi$ above, transforms as
$$M_k\act\varphi=0\quad,\quad N_k\act\varphi=mcdx_k\quad,\eqn\difxii$$
where $m$ is a mass parameter; thus, $\varphi$ has dimensions of an action (the
appearance of a new dimensional constant $m$ besides $c$ and $\kappa$, and
the form of
$N_k\act\varphi$ in \difxii, will be justified
later). Eqs. \difxii\ and \difix\ lead to
$$\eqalign{&N_k\act[x_i,\varphi]=-\delta_{ki}[x_0,\varphi]+mc[x_i,dx_k]\quad,
\cr &
N_k\act[x_0,\varphi]=-[x_k,\varphi]+mc[x_0,dx_k]-{mc\over\kappa}dx_k\quad,\cr
& M_k\act[x_i,\varphi]=\epsilon_{kil}[x_l,\varphi]\quad,\quad
M_k\act[x_0,\varphi]=0\quad.\cr}\eqn\difxiii$$
A solution to the corresponding system of equations with unknowns $A_{\mu
a}^b$ in \difv\ is given by (cf. \difxi )
$$\eqalign{&[x_0,dx_0]={1\over\kappa}dx_0+{1\over \kappa mc}\varphi\quad,\quad
[x_0,dx_i]=0\quad,\quad [x_0,\varphi]=-{\varphi\over\kappa}\quad,\cr
&[x_i,dx_0]={1\over\kappa}dx_i\quad,\quad [x_i,dx_j]=-\delta_{ij}{\varphi\over
\kappa m c}\quad,\quad [x_i,\varphi]=0\quad,\cr}\eqn\difxiv$$
which may be checked to satisfy eq. \difvi\ (Jacobi) and Leibniz's rule
(with $d\varphi=0$) for the last eqs. in \IIxx\ (eq. \difvii)
(in fact, the above
solution \difxiv\ is just a simple one in an existing one-parameter family of
solutions; this parameter is related to the scale invariance
$[\phi\to\alpha\phi]$ of $M_i\act\phi=0=N_i\act\phi$).
Then, if we now replace
$\kappa$ by $\tilde\kappa/c$ the nonrelativistic limit of
eqs. \difxi\ is determined  by $[x_i,x_j]=0
\;,\;[x_i,t]={x_i\over\tilde\kappa}$
plus the $c\rightarrow\infty$ limit of \difxiv, namely
$$\eqalign{&[t,dt]={dt\over\tilde\kappa}\ ,\ [t,dx_i]=0
\quad,\quad [t,\varphi]=
-{\varphi\over\tilde\kappa}\quad,\cr
&[x_i,dt]={dx_i\over\tilde\kappa}\quad,\quad
[x_i,dx_j]=-\delta_{ij}{\varphi\over
m\tilde\kappa}\quad,\quad[x_i,\varphi]=0\quad.\cr}\eqn\difxv$$
In the undeformed case (\kap,\ $\hat\kappa$
or \kapp\ $\rightarrow\infty$) all the expressions
become commuting ones, and $\varphi$ (or $\phi$ in \difxii) becomes
`uncoupled' to the spacetime variables. To see the meaning of $\varphi$, it is
convenient to look first at a larger, eleven generator,
\kap -deformed Poincar\'e Hopf algebra.

\chapter{Pseudoextended \kap-Poincar\'e algebra $\tilde{\cal P}_\kappa$}
We now consider a
deformation of the four-dimensional Poincar\'e Lie algebra
centrally `pseudoextended' by
a one-generator algebra. The word `pseudoextension' refers to the fact that,
although a Lie group $G$ may have trivial second cohomology group
$H^2(G,U(1))$,
one may describe the {\it direct product} extension by means of a
two-coboundary
which in the contraction limit leads
to a non-trivial two-cocycle of a {\it central} extension
$\tilde G_c$ of the contraction $G_c$ of $G$. Although the two-coboundary
is trivial, it is not completely so in the sense that in the contraction
it gives rise to a non-trivial cohomology group element: a {\it
contraction may generate group cohomology} \refmark{\AACON, \AI}.
This is the case for the four-dimensional
Poincar\'e group $P$, for which $P_c$ is the Galilei group and $\tilde G_c$
is the 11-parameter central extension of the Galilei group \refmark{\BAR},
usually denoted $\tilde G_{(m)}$ since for the Galilei group $H^2(G,U(1))=R$
and the mass parameter characterizes the extension two-cocycle. This
non-trivial two-cocycle is the contraction limit of a two-coboundary generated
by a one-cochain which does not have a
contraction limit \refmark{\AACON, \SAL}\refmark{\FOOTCOCH}.

In terms of the Lie
algebra generators, a pseudoextension may appear as the consequence of a
redefinition of the basis of the trivially extended algebra which involves the
contraction parameter. We accordingly introduce $
{\cal P}_\kappa\times{\cal U}(\Xi)$ as the Hopf algebra of generators
$(M_i, N_i, P_j, P_0^\prime, \Xi)$,
where $\Xi$ is central ($[\Xi,{\rm all}]=0,\;[\Xi]=({\rm action})^{-1}$),
defined by the Hopf algebra relations of \pk\ plus
$\Delta\Xi=\Xi\otimes 1+1\otimes\Xi,\;S(\Xi)=-\Xi$, $\epsilon(\Xi)=0$. If
we now make the change
$$P_0^\prime=P_0-mc\Xi\quad,\eqn\VIIi$$
the Hopf algebra structure is written as before with the replacement of
$P_0^\prime$
by $P_0-mc\Xi$. Explicitly, the deformed algebra
$\tilde{\cal P}_\kappa$ is defined by
 $$\eqalign{&[P_i,P_j]=0\quad,\quad [P_i,P_0]=0\quad ,\quad
[M_i,M_j]=\epsilon_{ijk}M_k\quad,\cr
&[M_i,P_j]=\epsilon_{ijk} P_k\quad,\quad
[M_i,P_0] =0\quad,\quad
[M_i,N_j] =\epsilon_{ijk} N_k\quad,\cr
&[N_i,P_0] =P_i\quad,\quad [N_i,N_j]=-\epsilon_{ijk}M_k\quad,\cr
&[N_i,P_j] =\delta_{ij}\left[{\kappa\over 2} \left(1-\exp\left(-2{P_0-mc\Xi
\over\kappa}\right)\right)+{1\over{2\kappa}}{\bf P}^2\right]-
{1\over\kappa} P_iP_j\quad; \cr}\eqn\VIIa$$
$$\eqalign{&
\Delta  P_i = P_i\otimes 1+\exp\left(-{P_0-mc\Xi\over
\kappa}\right)\otimes P_i\quad, \cr &
\Delta N_i = N_i\otimes 1+
\exp\left(-{P_0-mc\Xi\over\kappa}\right)
\otimes N_i +{1\over\kappa} \epsilon_{ijk}
 P_j\otimes M_k\quad; \cr}\eqn\VIIb$$
the other coproducts being primitive.
In the \kap$\rightarrow\infty$ limit, $[N_i,P_j]=\delta_{ij}(P_0-mc\Xi)$; if
the $c\rightarrow\infty$ limit is now taken we get the $\tilde{\cal G}_{(m)}$
commutators including the [boosts, momenta] commutator
$[V_i,P_j]=-m\delta_{ij}\Xi$ (these two limits cannot be
interchanged; there is no $c\rightarrow\infty$ limit for \VIIa, \VIIb\ due to
$mc\Xi$).
Let now $\chi$ be the coordinate dual to the central one $\Xi,\;<\Xi,\chi>=1$.
 Then eqs.
\IIxx\ are now completed with
$$\Delta\chi=\chi\otimes 1+1\otimes\chi\quad,\quad
[\chi,x_i]={mc\over\kappa}x_i\quad,\quad
[\chi,x_0]=0\quad\eqn\VIIii$$
and $\Delta_L(d\chi)=1\otimes\chi\;,\;\Delta_R(d\chi)=d\chi\otimes 1$.
If we now introduce an {\it enlarged
\kap- Minkowski spacetime} $\tilde{\cal
M}_\kappa$ with coordinates $(x^\mu,\chi),\;[\chi]=$action, the left
 action of the Lorentz generators on
$(x^\mu,\chi)$ is given by \IIxxi\ plus
$$M_i\act\chi=0\quad,\quad N_i\act\chi=mcx_i\quad,\eqn\VIIiii$$
which imply $M_i\act d\chi=0,\;N_i\act d\chi=mcdx_i$ (cf. \difxii);
clearly, $N_i\act(mcx_0+\chi)=0$. The action of the Lorentz generators on the
commutators involving differentials is now given by \difix\ plus \difxiii\ in
which the one-form $\varphi$ is replaced by $d\chi$ plus
$$\eqalign{&M_k\act[\chi,dx_i]=\epsilon_{kil}[\chi, dx_l]\quad , \quad
M_k\act[\chi,dx_0]=0\quad,\;M_k\act[\chi,d\chi]=0\quad,\cr
&N_k\act[\chi,dx_i]=mc[x_k,dx_i]-\delta_{ki}[\chi,dx_0]-{mc\over\kappa}
\delta_{ki}dx_0\cr
&N_k\act[\chi,dx_0]=mc[x_k,dx_0]-[\chi,dx_k]-{mc\over\kappa}dx_k\cr
&N_k\act[\chi,d\chi]=mc[x_k,d\chi]+mc[\chi,dx_k]+
{m^2c^2\over\kappa}dx_k\quad.\cr}
\eqn\VIIiv$$
It may now be seen that eqs. \difxiv\ with $\varphi=d\chi$ plus the commutators
in \VIIii\ and
$$[\chi,dx_0]=-{d\chi\over\kappa}\quad,\quad
[\chi,dx_i]={mc\over\kappa}dx_i\quad,\quad
[\chi,d\chi]={2mc\over\kappa}d\chi\eqn\VIIv$$
define a covariant first order differential calculus on $\tilde{\cal M}_\kappa$
which is a solution of the system of equations defined by \difix,\ \difxiii\
(with $\varphi=d\chi$),
and
\VIIiv, which satisfies \difvi\ and \difvii.
The origin of the additional one form $\varphi$ is now clear:
it is the differential
of the new variable $\chi$ in the enlarged spacetime.
As one might expect the redefinition
$\chi^\prime\equiv\chi+mcx_0$ ($N_i\act\chi^\prime=0$)
takes the solution \VIIii, \difxiv\ (with
$\varphi=d\chi$) to the solution \difxi\ where $\phi={1\over mc}
d\chi^\prime$ and $\chi^\prime$ is now scalar.
As for $m$, it characterizes the
coboundary implied by the redefinition \VIIi.
The differential calculi based on $\varphi$ and $\phi$ are equivalent in
Woronowicz's sense but they are inequivalent in the nonrelativistic limit
(see eqs. \difxv\ and sec. 5).

To conclude, let us show that the  pseudoextended
$\kappa$-Poincar\'e $\tilde{\cal
P}_\kappa$ Hopf algebra
has a cocycle bicrossproduct structure (see \refmark{\MAJSOB} and Appendix A),
where $H={\cal P}_\kappa$ and $A={\cal U}(\Xi)$.
In it, $\alpha$ and $\psi$ are taken to be trivial, $\alpha(a,h)=a\epsilon(h)$
and $\psi(h)=1\otimes 1\epsilon(h)$.
The mapping $\beta$ is defined by
$$\eqalign{\beta(P_0)=1\otimes P_0\quad,
&\quad\beta(P_i)=\exp(mc\Xi/\kappa)\otimes P_i\quad,\cr
\beta(M_i)=1\otimes M_i\quad,&\quad\beta(N_i)=\exp(mc\Xi/\kappa)\otimes
N_i\quad.\cr}\eqn\VIIva$$
and the coboundary by
$$\xi(N_i,P_j)=\delta_{ij}{\kappa\over 4}
\left(1-\exp\left({2mc\over\kappa}\Xi\right)\right)\quad.\eqn\VIIvi$$
Then, looking at in Appendix A it is seen that the consistency
formulae are satisfied and that the Hopf algebra structure of
$\tilde{\cal P}_\kappa$
is recovered; in particular, $[N_i, P_j]$ (eq. \VIIa )
and eqs. \VIIb\ are
recovered from \apxxix\ and \apxxx.

\chapter{A deformation $\tilde{\cal G}_{(m)\hat\kappa}$ of the extended Galilei
algebra $\tilde{\cal G}_{(m)}$}
A natural way of deriving a deformed Galilei Hopf algebra is to apply the
standard
$c\rightarrow\infty$ limit contraction \refmark{\IW} to \pk. Here, this means
 contracting with respect to the Hopf subalgebra defined by the translation
and rotation
generators using the
redefinitions
$$P_i =X_i\quad,\quad
N_i =cV_i\quad,\quad
M_i =J_i\quad,\quad
P_0 ={1\over c} X_t\quad,\eqn\rescal$$
the last one being solely motivated by the replacement of $x^0$ by $ct$ in
$P_0$.
The contraction is made with respect to a subalgebra and no new
deformation parameter appears. This is in contrast with the
contraction of ${\cal U}_q(so(2,3))$ leading to \pk,
which was done with respect to the Lorentz sector which is a subalgebra in
$so_q(2,3)$ only in the $q=1$ limit (the `boost' commutators in $so_q(2,3)$
give rise to momenta for $q\not = 1$) and hence
required that $q$ become involved in the contraction process by
setting $q=\exp(1/\kappa R)$.
If we perform this contraction in \poiv-\poiii, however, we obtain
the cocommutative, undeformed Hopf algebra
structure of the Galilei enveloping algebra ${\cal U}({\cal G})$.
Moreover, we already saw that the $c\rightarrow\infty$ limit cannot be taken
directly in \difa\ and in \difxi\ or \difxiv.
Thus, \kap\ {\it must} be accompanied by factors of $c$
if a non-relativistic limit has to be feasible.

Let us consider in this section the case a) in sec. 3, $\kappa\rightarrow
\hat\kappa
c\;,\; [\hat\kappa]=L^{-2}T$ and write $\tilde{\cal P}_{\hat\kappa}$ for
$\tilde{\cal P}_\kappa$ with $\kappa$ replaced by $\hat\kappa c$. Using
\rescal\ we obtain from eqs. \VIIa, \VIIb\ in the contraction limit the {\it
deformed extended Galilei Hopf algebra}
$\tilde{\cal G}_{(m)\hat\kappa}$ defined by
$$\eqalign{&[X_i,X_j]=0\quad,\quad[X_i,X_t]=0\quad,\quad
[J_i,J_j]=\epsilon_{ijk}J_k\quad,\cr
&[J_i,X_j]=\epsilon_{ijk}X_k\quad,\quad [J_i,X_t]=0\quad,\quad
[J_i,V_j]=\epsilon_{ijk}V_k\quad,\cr
&[V_i,X_t]=X_i\quad,\quad [V_i,X_j]=\delta_{ij}{\hat\kappa\over
2}(1-\exp(2m\Xi)/\hat\kappa)\quad,\quad [\Xi,\;{\rm all}]=0\quad,\cr}\eqn\Bi$$
$$\eqalign{&\Delta X_i=X_i\otimes 1+\exp(m\Xi/\hat\kappa)\otimes X_i\quad,\cr
&\Delta V_i=V_i\otimes 1+\exp(m\Xi/\hat\kappa)\otimes V_i\quad;\cr}\eqn\Bii$$
the other coproducts are primitive and the antipodes of the generators are
simply given by a change of sign but for
$$S(X_i)=-\exp(-m\Xi/\hat\kappa)X_i\quad,\quad
S(V_i)=-\exp(-m\Xi/\hat\kappa)V_i\quad.\eqn\antipod$$
The Casimir operators for $\tilde{\cal G}_{(m)\hat\kappa}$ are easily found.
They are
$$\eqalign{C_1=&X_t(1-\exp(2m\Xi/\hat\kappa))\hat\kappa-{\bf X}^2\quad,\cr
C_2=&\left[{\bf J}{\hat\kappa\over 2}\left(1-\exp(2m\Xi/\hat\kappa)\right)
-({\bf V}\times {\bf X})\right]^2
\quad,\quad C_3=\Xi\quad;\cr}\eqn\extras$$
In an undeformed quantum theory, we could set $\Xi\sim i/\hbar$ (the
dependence of the wavefunction on the central parameter may be factored out by
a $U(1)$-equivariance condition). Thus,
$C_1/2m\propto U_{\hat\kappa}$ and $C_2/m\propto{\bf S}^2_{\hat\kappa}$
constitute the
deformations of the internal energy and of the spin operators
(in nonrelativistic quantum mechanics, the position operator $\hat{\bf x}$\ is
equal to ${\bf V}/m$ once the factor $i\hbar$ is added).
The classical $\hat\kappa\rightarrow\infty$
limit of \Bi -\extras\ reproduces the Hopf algebra structure of the enveloping
algebra ${\cal U}(\tilde{\cal G}_{(m)})$ of the centrally extended Galilei
algebra $\tilde{\cal G}_{(m)}$ (in which $[V_i,X_j]=-m\delta_{ij}\Xi$),
and the standard $\tilde {\cal G}_{(m)}$ Casimir operators.

Let us now introduce a
differential calculus on the {\it enlarged $\hat\kappa$-Newtonian spacetime}
\nkk\ of coordinates $(t,x_i,\chi)$ associated with
$\tilde{\cal G}_{(m)\hat\kappa}$. The left actions
$$ J_i\act x_j=\epsilon_{ijk}x_k\quad\quad
J_i\act t=0\quad,\quad
V_i\act x_j=-\delta_{ij}t\quad,\quad
V_i\act t=0\quad,\eqn\Biii$$
plus
$$ J_i\act\chi=0\quad,\quad V_i\act\chi=mx_i\quad,\eqn\Biv$$
lead to
$$\eqalign{&V_k\act[x_i,dx_j]=-\delta_{ki}[t,dx_j]-\delta_{kj}[x_i,dt]\quad,
\quad V_k\act[t,dx_i]=-\delta_{ki}[t,dt]\quad,\cr
&V_k\act[x_i,dt]=-\delta_{ki}[t,dt]\quad,\quad
V_k\act [t,dt]=0\quad,\cr
&J_k\act[x_i,dx_j]=\epsilon_{kil}[x_l,dx_j]+\epsilon_{kjl}[x_i,dx_l]\quad,
\quad J_k\act[t,dx_i]=\epsilon_{kil}[t,dx_l]\quad,\cr
&J_k\act[x_i,dt]=\epsilon_{kil}[x_l,dt]\quad,\quad
J_k\act[t,dt]=0\quad;\cr}\eqn \Bv$$
$$\eqalign{&V_k\act [x_i,d\chi]=-\delta_{ki}[t,d\chi]+m[x_i,dx_k]\quad,\quad
V_k\act[t,d\chi]=m[t,dx_k]\quad,\cr
&J_k\act[x_i,d\chi]=\epsilon_{kil}[x_l,d\chi]\quad,\quad
J_k\act[t,d\chi]=0\quad;\cr}\eqn\Bvi$$
$$\eqalign{&V_k\act[\chi,dx_i]=m[x_k,dx_i]-\delta_{ki}[\chi,dt]-
{m\over\hat\kappa}\delta_{ki}dt\quad,\quad
V_k\act[\chi,dt]=m[x_k,dt]\quad,\cr
&V_k\act[\chi,d\chi]=m[x_k,d\chi]+m[\chi,dx_k]+{m^2\over\hat\kappa}dx_k\quad,
\cr & J_k\act[\chi,dx_i]=\epsilon_{kil}[\chi,dx_l]\quad,\quad
J_k\act[\chi,dt]=0\quad,\quad J_k\act[\chi,d\chi]=0\quad.\cr}\eqn\Bvii$$
It is not difficult to check that the commutators
$$\eqalign{&[x_i,x_j]=0\quad,\quad[t,x_i]=0\quad,\quad [t,\chi]=0\quad,\quad
[x_i,\chi]=-{m\over\hat\kappa}x_i\quad,\cr
&[x_i,dt]=0\quad,\quad [x_i,dx_j]=0\quad,\quad [t,dt]=0\quad,\quad
[t, dx_i]=0\quad,\cr
&[t,d\chi]=0\quad,\quad [x_i,d\chi]=0\quad,\cr
&[\chi,dx_i]={m\over\hat\kappa}dx_i\quad,\quad [\chi,dt]=0\quad,\quad
[\chi,d\chi]={2m\over\hat\kappa}d\chi\quad,\cr}\eqn\Bix$$
are a solution to \Biii, \Biv, \Bv, \Bvi\ and \Bvii\ which satisfies Leibniz
rule and the Jacobi identities. Since \Bv, \Bvi, \Bvii\ and \Bix\ are
the $c\rightarrow\infty$ limits of \difix, \difxiii, \VIIiv\ and of \IIxx,
\VIIii, \difxiv\ and \VIIv\ respectively, we see that there is complete
consistency among the contraction limit $\tilde{\cal
P}_{\hat\kappa}\rightarrow$\gkk\ and the
$c\rightarrow \infty$ limit  relating the differential calculi
$\Gamma(\tilde{\cal M}_{\hat\kappa})\;,\;
\Gamma(\tilde{\cal N}_{\hat\kappa})$ on the enlarged
$\hat\kappa$-Minkowski $\tilde{\cal M}_{\hat\kappa}$ and Newtonian $\tilde{\cal
N}_{\hat\kappa}$ spacetimes respectively
associated with $\tilde{\cal P}_{\hat\kappa}$
and $\tilde{\cal G}_{(m)\hat\kappa}$. Thus, the diagrams
$$
\vbox{\settabs \+1234567&1234567&123456789012&1234567&34567
&1234567&1234567&\cr
\+$\ \ \ \ \ \ \tilde{\cal P}_{\hat\kappa}$&
$c\to\infty\atop {\raise 5pt \hbox{$\longrightarrow$}}$&
\gkk&
$\ \Gamma(\tilde{\cal M}_{\hat\kappa})$&
$c\to\infty\atop {\raise 5pt \hbox{$\longrightarrow$}}$&
$\Gamma(\tilde{\cal N}_{\hat\kappa})$\cr
\+ ${\scriptstyle\hat\kappa\to\infty}\;\downarrow $&&
$\ \downarrow$&
${\scriptstyle\hat\kappa\to\infty}\;\downarrow$&&$\ \downarrow$\cr
\+$\ \ \ \ \ \ \tilde{\cal P}$&
$\atop {\raise 5pt \hbox{$\longrightarrow$}}$&
$\tilde{\cal G}_{(m)}$&
$\ \Gamma(\tilde{\cal M})$&$
\atop {\raise 5pt \hbox{$\longrightarrow$}}$&
$\Gamma({\tilde{\cal N}})$\cr}
$$
are commutative.
The deformation described by \gkk\ is rather mild: it only affects
$[V_i,X_j]$ and $\Delta X_i$, $\Delta V_i$ (eqs. \Bi, \Bii) and,
as far as the differential calculus is concerned, the commutators of $\chi$
with $x_i,\; dx_i$ and $d\chi$ only; in particular, time is commutative. This
is not surprising if one realizes that \gkk\ (eqs. \Bi, \Bii) provides an
example of a cocycle extended Hopf algebra, the non-trivial antisymmetric
two-cocycle (generated by the contraction process) being given by
$$\xi(V_i,X_j)=\delta_{ij}{{\hat\kappa}\over
4}(1-\exp({2m\Xi\over\hat\kappa}))=-\xi(X_j,V_i)\quad\eqn\Bx$$
\ie, by the $c\rightarrow\infty$ limit of \VIIvi; this reproduces $[V_i, X_j]$
in \Bi.

To complete the picture, we mention that we might have looked at the
nonrelativistic limit of ${\cal P}_{\hat\kappa}\times{\cal U}(\Xi)$ itself.
The differential calculus for it is given by \difxi\ completed with $[\eta,
x_\mu]=0\;,\;[\eta,dx_u]=dx_\mu/\kappa\;,\; [\eta,d\eta]=0$ where $
\phi=d\eta$.
Replacing \kap\ by ${\hat\kappa}c$, we see that in the
nonrelativistic limit $\eta$ `decouples' and
that all commutators are
zero but for \NUU. This nonrelativistic calculus on Newtonian spacetime is
thus noncommutative despite the fact that it is associated to the
`classical' Galilei algebra: in the $c\to \infty$ limit, ${\cal P}
_{\hat\kappa}\times{\cal U}(\Xi)\to {\cal U}({\cal G})\times{\cal U}(\Xi)$
(rather than $\tilde{\cal G}_{(m){\hat\kappa}}$)
since ${\hat\kappa}$ disappears in the contraction limit.
Nevertheless, it may be seen that the differential calculus based on ${\cal G}$
allows for a proportionality constant $\mu$ in $[x_i, dx_j]=\delta_{ij}\mu dt$
(all other commutators must be zero), and thus the commutativity of the
$c\rightarrow\infty$ limit is consistent with the above result
if $\mu$ is set equal to $1/\hat\kappa$ instead of being zero.

\vskip 24 pt
\chapter{Deformed Galilei algebra ${\cal G}_{\tilde\kappa}$
 and its bicrossproduct structure}
Let us now consider the redefinition b) in Sec. 3 \ie, the algebra ${\cal
P}_{\tilde\kappa}$ obtained from \pk\ by replacing \kap\ by \kapp$/c$.
It will turn out that it is not possible to construct fully commutative
diagrams as in the previous section. Nevertheless, the
$c\rightarrow\infty$ limit of ${\cal P}_{\tilde\kappa}$
(eqs. \poiv, \pov\ and \poiii) gives rise to the {\it deformed
\kapp-Galilei algebra} \gk\ \refmark{\LODZ}\refmark{\FOOTFACTOR}.
Since this is a deformation of the ten parameter Galilei algebra ${\cal G}$
which is interesting in itself we shall describe it now. Its
commutators are given by the abelian translation sector plus
$$\eqalign{&[J_i,J_j] =\epsilon_{ijk}J_k\quad,\quad
[J_i,X_j] =\epsilon_{ijk}X_k\quad,\quad
[J_i,X_t] =0\quad, \cr &
[J_i,V_j] =\epsilon_{ijk}V_k\quad,\quad
[V_i,X_t] =X_i\quad, \cr &
[V_i,X_j] =\delta_{ij}{1\over{2\tilde\kappa}}{\bf X}^2-{1\over\tilde\kappa}
X_iX_j\quad,\quad
[V_i,V_j] =0\quad; \cr &}\eqn\gali$$
the coproducts and antipodes are given by (cf. \pov, \poiii)
$$\Delta X_t =X_t\otimes 1+1\otimes X_t\quad,\quad
\Delta X_i =X_i\otimes 1 +\exp(-X_t/\tilde\kappa)\otimes X_i\quad,\eqn\galii$$
$$\Delta J_i =J_i\otimes 1+1\otimes J_i\quad,\quad
\Delta V_i =V_i\otimes 1+ \exp(-X_t/\tilde\kappa)\otimes V_i+{\epsilon_{ijk}
\over\tilde\kappa} X_j\otimes J_k\quad,\eqn\galiii$$
$$S(X_t) =-X_t\quad,\quad
S(X_i) =-\exp(X_t/\tilde\kappa) X_i\quad,\eqn\galiv$$
$$S(J_i) =-J_i \quad,\quad
S(V_i) =-\exp(X_t/\tilde\kappa) V_i +{1\over\tilde\kappa}\epsilon_{ijk}
\exp(X_t/\tilde\kappa) X_jJ_k\quad.\eqn\galv$$
Eqs. \gali -\galv\  satisfy all Hopf algebra
axioms; in the \kapp
$\rightarrow\infty$ limit, the undeformed Galilei
Lie algebra ${\cal G}$ expressions are
obtained.

Since ${\cal G}$ has a semidirect product structure with
the translations being an ideal, it is natural to ask
ourselves whether the \kapp-deformation above has a
bicrossproduct structure. We now show that this is the case, and that in fact
it
may be obtained by the contraction limits of the right action $\alpha$
and the left coaction $\beta$ in \aldef\ and \bedef. Specifically,
${\cal G}_{\tilde\kappa}=
{\cal U}(R\circ B)\bic {\cal U}_{\tilde\kappa}(Tr)$
where now
\mypointbegin
${\cal U}_{\tilde\kappa}(Tr)$ is the commutative and noncocommutative algebra
defined by \galii\  and \galiv
\mypoint
${\cal U}(R\circ B)$ is the undeformed Hopf algebra of
rotations and Galilean boosts,
$$[J_i,J_j]=\epsilon_{ijk}J_k\quad,\quad [J_i,V_j]=\epsilon_{ijk}V_k\quad,\quad
[V_i,V_j]=0\quad,\eqn\galvi$$
with antipodes $S(J_i)=-J_i,\  S(V_i)=-V_i$.
\mypoint
the right action
$\alpha:{\cal U}_\kappa(Tr)\otimes {\cal U}(R\circ B)\rightarrow
{\cal U}_\kappa(Tr)$
is defined by:
$$\eqalign{&X_t\acti J_i \equiv[X_t,J_i]=0\quad,\quad
X_t\acti V_i \equiv[X_t,V_i]=-X_i\quad, \cr &
X_i\acti J_j \equiv[X_i,J_j]=\epsilon_{ijk}J_k\quad,\cr &
X_i\acti V_j \equiv[X_i,V_j]=-\delta_{ij}{1\over 2\tilde\kappa}{\bf X}^2+
{1\over\tilde\kappa}X_iX_j\quad, \cr &}\eqn\defalfa$$
and the left coaction
$\beta:{\cal U}(R\circ B)\rightarrow
{\cal U}_{\tilde\kappa}(Tr)\otimes{\cal U}(R\circ B)$ by
$$\beta(J_i) =1\otimes J_i\quad,\quad
\beta(V_i) =\exp(-X_t/\tilde\kappa)\otimes V_i+
{\epsilon_{ijk}\over \tilde\kappa}X_j\otimes V_k\eqn\defbeta$$
\ie, by the $c\rightarrow\infty$ limits of \aldef,\ \bedef.

To prove that ${\cal G}_{\tilde\kappa}={\cal U}(R\circ B)\bic{\cal
U}_{\tilde\kappa}(Tr)$ one
needs checking first that the properties \apiii\ (module action)
\apv, \apvi\ (comodule coaction)
and \apvii,\ \apviii\ (comodule coalgebra) are satisfied with $\alpha$ and
$\beta$ defined by \defalfa, \defbeta\
and that the compatibility conditions \apix-\apxiii\ hold.
In the present case, the property \apiv\ for $\alpha$ is
 automatic since $\alpha$ is
given in terms of commutators, eq. \defalfa, and the coproduct in $R\circ B$ is
primitive. Moreover, it is sufficient to
check \apvii,\ \apviii\ for the elements $h=V_i\in {\cal U}(R\circ B)$ since
for $\beta$ trivial ($\beta(h)=1_A\otimes h$ for $h=J_i$)
and primitive coproducts ($\Delta(h)=h\otimes 1 +1\otimes h$) eqs. \apvii\
and \apviii\ are automatically satisfied.
As for the compatibility conditions, eqs. \apix\ and \apxi\ are
immediate, and \apxiii\ is automatic since $A\ (Tr)$ is abelian and $H$
(rotations and boosts) cocommutative (see \eg\ \refmark{\MB} or Appendix A).
As for the r.h.s. of
\apx, it reads
for \eg, $a=X_i,\ h=V_j$
$$[X_i,V_j]\otimes 1+[\exp(-X_t/\tilde\kappa),V_j]\otimes X_i+
\exp(-2X_t/\tilde\kappa)\otimes[X_i,V_j]+$$
$$+\exp(-X_t/\tilde\kappa){\epsilon_{jkl}\over\tilde\kappa}
X_k\otimes [X_i,J_l]\quad,\eqn\demiii$$
where $a\acti 1=a$ [\apiii]\ and $1\acti h=\epsilon(h)=0$ [\apiv]\
have been used and
which may be checked to be equal to $\Delta([X_i,V_j])$ as computed
from the last equation in \defalfa. Similarly, condition \apxii\ is satisfied
by \defbeta. The Hopf structure now follows from \apxiv-\apxvii.
Clearly, these equations
 do not modify the coproduct and antipode of the translation
sector and reproduce trivially those for $J_i$ since, for it, $\beta$ is
trivial. As for $V_i$, it is simple to check that eqs. \galiii\ and
\galv\ are obtained.

Finally, let us find the Casimir operators for \gk.
The Casimirs of \pk\ are given by
$$\eqalign{&C_1={\bf P}^2\exp(P_0/\kappa)-4\kappa^2\sinh^2(P_0/\kappa)\quad,\cr
&C_2=\left(\cosh(P_0/\kappa)-{{\bf P}^2\exp(P_0/\kappa)\over 4\kappa^2}\right)
W_0^2-{\bf W}^2\quad,\cr}\eqn\casi$$
where the deformed Pauli-Luba\'nski vector is given by
$$\eqalign{&W_0\equiv{\bf PM}\exp(P_0/2\kappa)\quad,\cr &
 W_i\equiv \kappa M_i\sinh(P_0/
\kappa)+
\exp(P_0/\kappa)\left(\epsilon_{ijk}P_jN_k+{1\over 2\kappa}(M_i{\bf
P}^2-P_i(\bf{PM}))\right)\quad.\cr}\eqn\casii$$
Taking the $c\rightarrow\infty$
limit of \casi,\ \casii, the Casimirs of \gk\ are
found to be
$$C_1={\bf X}^2\exp(X_t/\tilde\kappa)\quad,\quad
C_2=-
\exp(2X_t/\tilde\kappa)\left[{\bf X}\and{\bf V}
+{1\over 2\tilde\kappa}{\bf J\;X}^2\right]^2\quad,\eqn\casiii$$
where for the second one we have used $C^\prime\equiv C_2/c^2$ to take the
limit; it may be checked that they commute
with all elements of \gk.
Moreover, for $\tilde\kappa\rightarrow\infty$ they
 become the Galilei algebra Casimirs \ie, the square of momentum ${\bf
X}^2$ and of angular momentum $({\bf V}\and {\bf X})^2\propto {\bf L}^2$
respectively.

\chapter{\kapp-Newtonian spacetime ${\cal N}_{\tilde\kappa}$
and differential calculus}
The {\it \kapp-Newtonian spacetime} ${\cal N}_{\tilde\kappa}$
 may be introduced (as ${\cal M}_\kappa$) by duality, now from
 \galii. This leads to the basic relations
$$\eqalign{&\Delta t=t\otimes 1+1\otimes t\quad,\quad \Delta x_i=x_i\otimes 1+
1\otimes x_i\quad;
\cr &[x_i,x_j]=0\quad,\quad [x_i,t]={x_i\over \tilde\kappa}\quad,
\cr}\eqn\vi$$
analogous to \IIxx; again we find for this case, as for ${\cal M}_\kappa$,
the physically rather inconvenient fact of having a noncommutative
time. The left action of the \kapp-Galilei algebra generators is given by
\Biii.
As one might expect, the elements $x_0$ and
${\bf x}^2$ are invariant in ${\cal N}_{\tilde\kappa}$ under the
\kapp-Euclidean Hopf algebra generated by $(X_t,X_i,J_l)$, eqs. \gali-\galv.
To introduce a first order
 \kapp-Galilei invariant differential calculus
we apply \Biii, \difviii\ and obtain
$$\eqalign{&V_k\act [x_i,dx_j]=-\delta_{ki}[t,dx_j]-\delta_{kj}[x_i,t]+
{1\over\tilde\kappa}(\delta_{kj}dx_i-\delta_{ij}dx_k)\quad,\cr
&V_k\act[t,dx_i]=-\delta_{ki}[t,dt]+{1\over\tilde\kappa}\delta_{ki}dt\quad,\cr
&V_k\act[x_i,dt]=-\delta_{ki}[t,dt]\quad,\quad V_k\act[t,dt]=0\quad,\cr
&J_k\act[x_i,dx_j]=\epsilon_{kil}[x_l,dx_j]+\epsilon_{kjl}[x_i,dx_l]\quad,\cr
&J_k\act[t,dx_i]=\epsilon_{kil}[t,dx_l]\quad,\quad
J_k\act[x_i,dt]=\epsilon_{kil}
[x_l,dt]\quad,\quad J_k\act[t,dt]=0\cr}\eqn\viii$$
\ie, the contraction limit of \difix.

If we now try to find an expression for the commutators
 $[t,dt],\quad
[t,dx_i],\break\quad[x_i,dt],\quad[x_i,dx_j]$
following the same process which lead to
\difxi\ we find that there is no solution even if
an invariant one-form $\phi$ is added. Let us then introduce a one-form
$\varphi$, with dimensions of an action, and with transformation properties
$$V_k\act\varphi=mdx_k\quad,\quad J_k\act\varphi=0\quad.\eqn\viv$$
Using \viv, the relations \viii\ are now completed with the following ones
$$\eqalign{&V_k\act[x_i,\varphi]=-\delta_{ki}[t,\varphi]+m[x_i,dx_k]\quad,\cr
&V_k\act[t,\varphi]=m[t,dx_k]-{m\over\tilde\kappa}dx_k\quad,\cr
&M_k\act[x_i,\varphi]=\epsilon_{kil}[x_l,\varphi]\quad,\quad M_k\act
[t,\varphi]=0\quad,\cr}\eqn\vv$$
which turn out to be the contraction limit of eqs. \difxiii.
Writing \eg\  $[x_i,dx_j]=
A_{ij}^a\chi_a$ where $\chi_a=(dt,dx_i,\varphi)$, eqs. \viii\ and \vv\ give
rise to a linear system of equations, which admits \difxv\ as a solution.
The general solution depends on one parameter $\lambda$,
$$\eqalign{&[t,dt]={1\over\tilde\kappa}dt\quad,\quad [t,dx_i]=0\quad,\quad
[t,\varphi]=2m\lambda dt-{1\over\tilde\kappa}\varphi\quad,\cr
&[x_i,dt]={1\over\tilde\kappa}dx_i\quad,\quad [x_i,dx_j]=\delta_{ij}
\left(\lambda dt-{\varphi\over m\tilde\kappa}\right)\quad,\quad
[x_i,\varphi]=\lambda mdx_i\quad;}\eqn\vvi$$
eqs. \difxv\ correspond to $\lambda=0$. It is reasonable to select
$\lambda=0$ since $\lambda$ has dimensions, $[\lambda]=L^2T^{-1}$, and
there are no grounds to introduce another dimensionful parameter.
We do not have complete closure in this case however, because the last two
equations in \VIIv\ do not have a limit if \kap\ is replaced by \kapp$/c$ in
them.

\chapter{Conclusions and outlook}
We have given in this paper a deformation  \gkk\
  of the extended Galilei
algebra  $\tilde{\cal G}_{(m)}$  and a deformation \gk\
of the Galilei algebra ${\cal G}$, and discussed their
differential calculus on the enlarged Newtonian spacetime  $\tilde{\cal
N}_{\hat\kappa}$ and on ${\cal N}_{\tilde\kappa}$. The two deformations
have been
obtained, respectively, as the contraction limits of a pseudoextension
$\tilde{\cal P}_{\hat\kappa}$ of ${\cal P}_{\hat\kappa}$
and of ${\cal P}_{\tilde\kappa}$ .
Both \gkk\ and \gk\ retain the same cocycle and bicrossproduct
structure of their parent deformed algebras. In the case of $\tilde{\cal
P}_{\hat\kappa}$  and \gkk,
there is complete commutativity among the nonrelativistic and the
undeformed $(\hat\kappa\rightarrow\infty)$ limits.
The fact that an additional variable $\chi$ ($\varphi=d\chi$) is
necessary in the relativistic differential calculus associated with
$\tilde{\cal P}_{\hat\kappa}$,
and that the deformation
enters in \gkk\ only through the central generator, opens an intriguing
relation among deformation and quantization. Indeed, it is known that in
the undeformed case the central additional generator plays a r\^ole in
geometric quantization theories in which there exists a $U(1)$-principal
bundle structure. In them,
Planck's constant appears as the divisor which makes of the two-cocycle (local
exponent) the dimensionless quantity needed for a phase
(in particular, this would also be the situation for the simplest
example of the Weyl-Heisenberg algebra).
In any case, the need for this additional one-form in the presence of a
deformation is not an isolated fact; for instance, it
is also present in the case of the Euclidean space obtained starting from
${\cal U}_\omega{\cal E}(2)$, for which there is also a
bicrossproduct structure and a similar study can be made.
This phenomenon is similar to the known
unbalance between the invariant vector fields and Maurer-Cartan
one-forms present in deformed groups
other than the general linear groups. In fact, it appears to be difficult to
construct noncommutative differential `spaces' whose tangent spaces have the
same dimension as in the undeformed theory (see, \eg\ \refmark{\CWSWW},
\refmark{\ZSAL}, \refmark{\SCH} and \refmark{\AS}).

The commutativity of the $c\rightarrow\infty$ and undeformed
 limits fails for the
$\tilde{\cal P}_{\tilde\kappa}$ covariant differential calculus; this is also
manifest in the absence of an extended
$\tilde{\cal G}_{(m)\tilde\kappa}$-type deformation
coming from a $c\rightarrow\infty$ limit of $\tilde{\cal P}_{\tilde\kappa}$.
Thus, although a deformation of the Galilei algebra \gk\ and a deformed
differential calculus on $\tilde\kappa$-Newtonian spacetime ${\cal
N}_{\tilde\kappa}$ exist if the one-form $\varphi$ [\viv] is added,
it is not possible to define an enlarged
nonrelativistic spacetime $\tilde{\cal N}_{\tilde\kappa}$, since the presence
of $\chi$ does no allow for a $c\rightarrow\infty$ limit (see \VIIv).
The bicrossproduct Hopf algebra \gk\ is a stronger deformation
of the Galilei algebra (the time here is non-commutative, eq. \vi) than \gkk\
is of $\tilde{\cal G}_{(m)}$, but no cocycle bicrossproduct extension seems
to exist for \gk. Thus,
although the need for an additional form appears natural due to the central
generator in $\tilde{\cal P}_{\hat\kappa}$ and \gkk, no such extension exists
for \gk\ closing the appropriate commutative diagrams.

\ack
The authors wish to thank A. Ballesteros, P. P. Kulish, J. Lukierski, M. del
Olmo, M. Santander and F. Herranz for helpful discussions or comments.
This paper has been partially supported by a CICYT (Spain) research grant.
One of the authors (JCPB) wishes to thank the Spanish Ministry of
Education and Science and the CSIC for a grant.

\endpage

\APPENDIX{A}{A: Bicrossproduct of Hopf algebras and cocycles}
We summarize here for completeness Majid's bicrossproduct construction
and refer to \refmark{\MB,\MAJSOB} (see also \refmark{\BCM}) for details.
Let $A$ and $H$ be Hopf algebras, and let
\mypointbegin
$A$ be a right $H$-module algebra ($H\rimo A$)
\mypoint
$H$ be a left $A$-comodule coalgebra ($H\leco A$)
\ie, there exist linear mappings
$$\alpha:A\otimes H\rightarrow A\quad,\quad\alpha(a\otimes h)\equiv a\acti h
\quad,\quad a\in A,\ h\in H\quad;\eqn\api$$
$$\beta:H\rightarrow A\otimes H\quad,\quad \beta(h)=h^{(1)}\otimes
h^{(2)}\quad,
\quad h^{(1)}\in A,\ h^{(2)}\in H\eqn\apii$$

\noindent
(in general, superindices refer to $\beta$,
subindices to the coproduct $\Delta$), such that the properties of

\noindent
a1) $\alpha$ being a right $H-$module action:
$$a\acti 1_H=a\quad,\eqn\apiiia$$
$$(a\acti h)\acti h^{\prime}=a\acti hh^{\prime}\quad;
\eqn\apiii$$
a2) $A$ being a right H-module algebra:
$$1_A\acti h=1_A\epsilon(h)\quad,\quad
(ab)\acti h=(a\acti h_{(1)})(b\acti h_{(2)})\quad;
\eqn\apiv$$
b1) $\beta$ being a left A-comodule coaction:
$$\epsilon_A(h^{(1)})\otimes h^{(2)}=1_A\otimes h\equiv h\quad
[(\epsilon\otimes id)\circ\beta=id]\quad,\eqn\apv$$
$$h^{(1)}\otimes h^{(2)(1)}\otimes h^{(2)(2)}=h^{(1)}_{\ (1)}
\otimes h^{(1)}_{\ (2)}\otimes h^{(2)}\quad [(id\otimes\beta)\circ\beta=
(\Delta\otimes id)\circ\beta]\quad;\eqn\apvi$$
b2) $H$ being a left A-comodule coalgebra:
$$h^{(1)}\epsilon_H(h^{(2)})=1_A\epsilon_H(h)\quad
[(id\otimes \epsilon)\circ \beta=\epsilon]\quad,\eqn\apvii$$
$$\eqalign{&
h^{(1)}\otimes h^{(2)}_{\ (1)}\otimes h^{(2)}_{\ (2)}=
h_{(1)}^{\ (1)}h_{(2)}^{\ (1)}\otimes h_{(1)}^{\ (2)}
\otimes h_{(2)}^{\ (2)}\quad\cr [(id\otimes\Delta)\circ\beta&=(m_A\otimes id
\otimes id)\circ(id\otimes\tau\otimes id)\circ(\beta\otimes\beta)\circ\Delta
\equiv(\beta\hat\otimes\beta)\circ\Delta]\;,\cr}\eqn\apviii$$
where $m_A$ is the multiplication in $A$ and $\tau$ is the twist mapping,
are fulfilled.

Then, if the compatibility conditions
$$\epsilon_A(a\acti h)=\epsilon_A(a)\epsilon_H(h)\quad,\eqn\apix$$
$$\Delta(a\acti h)\equiv (a\acti h)_{(1)}\otimes (a\acti h)_{(2)}
=(a_{(1)}\acti h_{(1)})h_{(2)}^{\ (1)}\otimes a_{(2)}\acti
h_{(2)}^{\ (2)}\quad,\eqn\apx$$
$$\beta(1_H)\equiv 1_H^{(1)}\otimes 1_H^{(2)}=1_A\otimes 1_H\quad,\eqn\apxi$$
$$\beta(hg)\equiv (hg)^{(1)}\otimes (hg)^{(2)}
=(h^{(1)}\acti g_{(1)})g_{(2)}^{\ (1)}\otimes h^{(2)}g_{(2)}^{\ (2)}
\quad,\eqn\apxii$$
$$h_{(1)}^{\ (1)}(a\acti h_{(2)})\otimes h_{(1)}^{\ (2)}=(a\acti h_{(1)})
h_{(2)}^{\ (1)}\otimes h_{(2)}^{\ (2)}\quad,\eqn\apxiii$$
are satisfied \refmark{\FOOTCOND}, there is a Hopf algebra structure on
$K=H\otimes A$ called the (right-left) bicrossproduct $H\bic A$ \refmark{\MB}
defined by
$$(h\otimes a)(g\otimes b)=hg_{(1)}\otimes (a\acti g_{(2)})b\quad,
\quad h,g\in H
\;a,b\in A\quad,\eqn\apxiv$$
$$\Delta_K(h\otimes a)=h_{(1)}\otimes h_{(2)}^{\ (1)}a_{(1)}
\otimes h_{(2)}^{\ (2)}\otimes a_{(2)}\quad,\eqn\apxv$$
$$\epsilon_K=\epsilon_H\otimes\epsilon_A\quad,\quad
 1_K=1_H\otimes 1_A\quad,\eqn\apxvi$$
$$S(h\otimes a)=(1_H\otimes S_A(h^{(1)}a))(S_H(h^{(2)})\otimes 1_A)
\quad.\eqn\apxvii$$
In $H\otimes A,\;
h\equiv h\otimes 1_A$ and $a\equiv 1_H\otimes a$; thus,
$ah=h_{(1)}\otimes(a\acti h_{(2)})$.
When $\beta=1_A\otimes I$ \ie\
$\beta(h)=1_A\otimes h$ (trivial coaction) and $H$ is cocommutative, $K$ is
the semidirect {\it product}
 of Hopf algebras; when $\alpha$ is trivial, $\alpha=
1_A\otimes\epsilon_H\ (a\acti h=a\epsilon_H(h))$ and $A$ is commutative
$K$ is the semidirect {\it coproduct} of Hopf algebras\refmark{\MOL,\MB}.
When $\alpha$ is trivial, $\beta(hg)=\beta(h)\beta(g)$ (algebra homomorphism)
since $\epsilon(g_{(1)})\beta(g_{(2)})=\beta(g)$ (linearity of $\beta$).

The above construction may be now extended to include cocycles
\refmark{\MAJSOB}. Let $H$ and $A$ two Hopf algebras and $\alpha$ and $\beta$
  as in \api, \apii. Then $A$
is a right $H$-module cocycle algebra if \apiiia, \apiv\ are fulfilled and
there is a linear (two-cocycle) map $\xi:H\otimes H\rightarrow A$ such that
$$\xi(h\otimes 1_H)=1_A\epsilon(h)=\xi(1_H\otimes h)\quad [\xi(1_H\otimes
1_H)=1_A]\quad,\eqn\apxx$$
$$\xi(h_{(1)}g_{(1)}\otimes f_{(1)})
(\xi(h_{(2)}\otimes g_{(2)})\acti f_{(2)})
=\xi(h\otimes g_{(1)}f_{(1)})\xi(g_{(2)}\otimes f_{(2)})\;,\; \forall
h,g,f\in H\;,\eqn\apxxi$$
(cocycle condition \refmark{\FOOTASSOC})
and \apiii\ is replaced by
$$\xi(h_{(1)}\otimes g_{(1)})((a\acti h_{(2)})\acti
g_{(2)})=(a\acti(h_{(1)}g_{(1)}))\xi(h_{(2)}\otimes
g_{(2)})\quad, \forall a\in
A,\forall h,g\in H\quad,\eqn \apxxii$$
which for $\xi$ trivial reproduces \apiii. Similarly, $H$ is a left
$A$-comodule coalgebra cocycle if \apv,\ \apvii,\ \apviii\ are fulfilled, and
there is a linear map $\psi:H\rightarrow A\otimes A$,
$\psi(h)=\psi(h)^{(1)}\otimes\psi(h)^{(2)}$, such that
$$\epsilon(\psi(h)^{(1)})\psi(h)^{(2)}=1\epsilon(h)=\psi(h)^{(1)}\epsilon(\psi
(h)^{(2)})\ ,\ \left[(\epsilon\otimes
id)\circ\psi=(id\otimes\epsilon)\circ\psi\right]\ ,\eqn\apxxiii$$
$$h_{(1)}^{\ (1)}\psi(h_{(2)})^{(1)}\otimes\psi(h_{(1)}^{\ (2)})\Delta \psi
(h_{(2)})^{(2)}=\psi(h_{(1)})\Delta\psi(h_{(2)})^{(1)}\otimes
\psi(h_{(2)})^{(2)},\ \forall h\in H\quad,\eqn\apxxiv$$
(dual cocycle condition) and \apvi\ is replaced by
$$\eqalign{((id\otimes\beta)\circ\beta(h_{(1)}))(\psi(h_{(2)})\otimes 1)&=
(\psi(h_{(1)})\otimes 1)((\Delta\otimes id)\circ\beta(h_{(2)}))=\cr
&\psi(h_{(1)})
\Delta h_{(2)}^{\ (1)}\otimes h_{(2)}^{\ (2)}\quad.\cr}\eqn\apxxv$$

Then, if the compatibility conditions \apix,\ \apxi,\ \apxiii\ and
$$\psi(h_{(1)})\Delta(a\acti h_{(2)})=[(a_{(1)}\acti h_{(1)})h_{(2)}^{\ (1)}
\otimes a_{(2)}\acti h_{(2)}^{\ (2)}]\psi(h_{(3)})\quad, \eqn\apxxvi$$
$$\beta(h_{(1)}g_{(1)})(\xi(h_{(2)}\otimes g_{(2)})\otimes 1)=
\xi(h_{(1)}\otimes g_{(1)})
(h_{(2)}^{\ (1)}\acti g_{(2)})g_{(3)}^{\ (1)}\otimes h_{(2)}^{\
(2)}g_{(3)}^{\ (2)}\quad,\eqn\apxxvii$$
(which replace \apx,\apxii), together with
$$\eqalign{\psi(h_{(1)}g_{(1)})\Delta\xi(h_{(2)}\otimes g_{(2)})
&=\left[\xi(h_{(1)}\otimes g_{(1)})
(h_{(2)}^{\ (1)}\acti g_{(2)})g_{(3)}^{\ (1)}(\psi(h_{(3)})^{(1)}\acti g_{(4)})
g_{(5)}^{\ (1)}\right.\cr
&\left.\otimes\xi(h_{(2)}^{\ (2)}\otimes g_{(3)}^{\ (2)})(\psi(h_{(3)})^{(2)}
\acti g_{(5)}^{\
(2)})\right]\psi(g_{(6)})\quad,\cr}\eqn\apxxviii$$
$$\eqalign{\epsilon(\xi(h\otimes g))&=\epsilon(h)\epsilon(g)
\quad,\quad \psi(1_H)=1_A\otimes 1_A\cr}\eqn\apxxviiia$$
hold, $(H,A,\alpha,\beta,\xi,\psi)$ determine a cocycle right-left
bicrossproduct bialgebra \break
$H ^\psi\bic_\xi A$ \refmark{\MAJSOB}. In it, the counit and
unit are defined by \apxvi\ and the product and coproduct by
$$(h\otimes a)(g\otimes b)=h_{(1)}g_{(1)}\otimes \xi(h_{(2)}\otimes
g_{(2)})(a\acti
 g_{(3)})b \quad,\eqn\apxxix$$
$$\Delta(h\otimes a)=h_{(1)}\otimes h_{(2)}^{\ (1)}\psi(h_{(3)})^{(1)}a_{(1)}
\otimes h_{(2)}^{\ (2)}\otimes \psi(h_{(3)})^{(2)}a_{(2)}\quad.\eqn\apxxx$$

For $\xi$ trivial $[\xi(h\otimes g)=\epsilon(h)\epsilon(g)1_{A}]$
\apxxii\ reduces
to \apiii, \apxxvii\ to \apxii\ (use $(\epsilon\otimes id\otimes
id)\Delta^2=\Delta)$ and \apxxix\ to \apxiv.
For $\psi$ trivial $[\psi(h)=1_{A}\otimes 1_{A}\epsilon(h)]$, \apxxv\
reduces to \apvi, \apxxvi\ to \apx\ (use $(m\otimes id)(id\otimes \beta)
(id\otimes id\otimes\epsilon)\Delta^2=(m\otimes id)(id\otimes\beta)\Delta)$
and \apxxx\ to \apxv\ (use $(id\otimes
 id\otimes\epsilon)\Delta^2=\Delta$ multiplied from the right by
$(id\otimes\Delta a)$).

\endpage
\refout
\endpage
\end